\documentclass%[10pt]
{amsart}
\usepackage{amsmath,amsthm,amsfonts,amscd,amssymb,latexsym}
\usepackage{graphicx}

\def\sm{\setminus}
\def\theta{\vartheta}
\def\disk{{\mathbb D}}

\def\ovl{\overline}
\def\phi1{\phi}
\def\phi{\varphi}
\def\eps{\varepsilon}
\def\theta{\vartheta}

\def\Im{\mbox{\rm Im}}

\def\sm{\setminus}

\def\N{\mathbb {N}}

\def\Z{\mathbb {Z}}
\def\R{\mathbb {R}}
\def\C{\mathbb {C}}
\def\Cbar{\overline{\C}}

\def\disk{\mathbb {D}}

\def\diskbar{\ovl{\disk}}

\newtheorem{theorem}{Theorem} 

\newtheorem{proposition}[theorem]{Proposition}
\newtheorem{lemma}[theorem]{Lemma}
\newtheorem{definition}[theorem]{Definition}

\newtheorem{corollary}[theorem]{Corollary}

\theoremstyle{definition}
\newtheorem*{remark}{Remark}

\newcommand{\Sd}{{\mathcal S}_d}
\newcommand{\Pd}{{\mathcal P}_d}

\def\lineclear
    {\rule{0pt}{0pt}\nopagebreak\par\nopagebreak\noindent}

\def\reminder #1 {{\sf #1}}
\def\hide #1 {}

\newcommand{\comment}[1]{\marginpar{#1}}
\renewcommand{\comment}[1]{}

\title[Efficient GLobal Newton Dynamics]{On the Efficient Global Dynamics of Newton's Method for Complex Polynomials}

\author{Dierk Schleicher}
\address{Jacobs University Bremen, Postfach 750 561, D-28725 Bremen, Germany}
\email{dierk@jacobs-university.de}

\begin{document}

\begin{abstract}
We investigate Newton's method as a root finder for complex polynomials of arbitrary degree. While polynomial root finding continues to be one of the fundamental tasks of computing, with essential use in all areas of theoretical mathematics, numerics, computer graphics and physics, known methods have either excellent theoretical complexity but cannot be used in practice, or are practically efficient but are a lacking a successful theory behind them. 

In this manuscript we investigate the theoretical complexity of Newton's method for finding all roots of polynomials of given degree and show that it is near-optimal for the known set of starting points that find all roots. This theoretical result is complemented by a recent implementation of Newton's method that finds all roots of various polynomials of degree more than a million, significantly faster than our upper bounds on the complexity indicate, and often much faster than established fast root finders. In some experiments, it was possible to find all roots using Newton's method even with complexity $O(d\log d)$ for degrees exceeding 100 million. Newton's method thus stands out as a method that has merits both from the theoretical and from the practical point of view.

Our study is based on the known explicit set of universal starting points, for each degree $d$, that are guaranteed to find all roots of polynomials of degree $d$ (appropriately normalized). We show that this set contains $d$ points that converge very quickly to the $d$ roots: the expected total number of Newton iterations required to find all $d$ roots with precision $\eps$ is $O(d^3\log^3d+d\log|\log\eps|)$, which can be further improved to $O(d^2\log^4d+d\log|\log\eps|)$; in the worst case allowing near-multiple roots, the complexity is $O(d^4\log^2d+d^3\log^2d|\log\eps|)$.   The arithmetic complexity for all these estimates is the same as the number of Newton iterations steps, up to a factor of $\log^2 d$.
\end{abstract}

\maketitle

\section{Introduction}

Finding roots of polynomials is one of the oldest problems in mathematics, and it is of significant interest today, in all areas of theoretical mathematics as well as in applications such as computer algebra and computer geometry; especially in statistical physics and dynamical systems, there is a natural need to find all roots of polynomials of very high degrees.
Surprisingly, theory and practice are not as well understood as one might expect.

There are root-finding algorithms with near-optimal theoretical complexity, notably by Pan \cite{Pan_NearOptimal}, but with constants so big that they cannot be used in practice. Then there are known and established practical implementations, notably \texttt{MPSolve 3.0} by Bini and Robol \cite{BiniRobol}, based on iteration in several variables (the Aberth-Ehrlich method), but they are lacking theory and have neither a proof of convergence in general nor an estimate on their speed of convergence --- but they work well in practice. Eigenvalue methods perform well (see for instance \cite{BiniEigenvalue}) especially for moderate degrees.

Newton's root-finding method is as old as analysis, and it is known to be a very efficient method for \emph{locally} finding approximate roots of smooth equations such as polynomials: once a reasonable approximation to a simple root is known, every iteration of the Newton method doubles the number of valid digits. However, Newton's method has a reputation as being difficult to understand as a \emph{global} dynamical system due to its ``chaotic'' nature.

In this paper we lay the foundations for establishing Newton's method as a root finder that is supported by theory and that works well in practice. It is known from \cite{NewtonInventiones} that for every degree $d$ there is a small universal set $\Sd$ of starting points that finds all roots of all complex polynomials of degree $d$ (appropriately normalized) when Newton's method is started at these points; we have $|\Sd|=1.11d\log^2d$. In this paper (with a later refinement in \cite{NewtonTodor}), we show that there is a subset of $d$ of these points that together only need $O(d^2\log^4d+d\log|\log\eps|)$ Newton iterations to find all roots with precision $\eps$ at least in the absence of near-multiple roots (if there are near-multiple roots, then all the isolated roots are found with this speed; the others are found too but with complexity $O(d^3\log^2d(d+|\log\eps|)$\,). These complexity bounds are near-optimal (up to logarithmic factors in $d$) for methods that start the Newton iteration at uniform distance from the disk containing all roots. 

\looseness-2
Our theoretical estimates are complemented by practical experiments~\cite{NewtonRobin,NewtonRobin2} that show that Newton's method routinely finds all roots of complex polynomials of degrees up to 134 million, and under certain conditions even significantly faster than the established root finder \texttt{MPSolve 3.0}. 

Newton's method as a global root finder of polynomials is an iterated rational map and thus indeed ``chaotic'' on its Julia set; however, methods from holomorphic dynamics may be brought to bear to control the dynamics. Here are some of the challenges that Newton's method faces even in the case of a polynomial $p$ in a single complex variable:
\begin{itemize}
\item
orbits of the Newton map that get close to zeroes of the derivative $p'$ will, under the Newton dynamics, jump near $\infty$ and will take a long time until they can get close to roots (if ever);
\item
there may be open sets in $\C$ in which the Newton dynamics does not converge to any root of $p$ (even for as simple polynomials as $p(z)=z^3-2z+2$);
\item
the boundary of the attracting basins of the various roots may have positive measure, so that the set of ``bad'' starting points has positive measure;
\item
even if almost all starting points in $\C$ converge to \emph{some} root of $p$, it is not clear to find starting points for all roots: it is conceivable that some roots can be found only from small sets of starting points (we do not want to use deflation: it is numerically unstable unless the roots are found in a certain order, and it may destroy specific forms of the polynomial that are easy to evaluate);
\item
finally, even if all roots are found, one needs efficient estimates on the required number of iterations.
\end{itemize}
This manuscript addresses all these issues: we specify, for arbitrary polynomials in a single complex variable, a universal set of starting points (depending only on the degree and some normalization) from which \emph{all} roots of \emph{all} polynomials of given degree are found, and so that the required number of iterations (or arithmetic complexity) is small: it is $O(d^4)$ in the worst case, and $O(d^3)$ or even $O(d^2)$ on average (up to factors of $\log d$). More precisely, we will prove the following theorem (the first half of which is not new, but required to state the main result):

\begin{theorem}[Newton efficiency]
\lineclear
For every degree $d\ge 2$, let $\Pd$ be the set of complex polynomials of degree $d$, normalized so that all roots are in the complex unit disk $\disk$. Then there is an explicit and finite universal set $\Sd$ consisting of $3.33\,d\log^2d(1+o(1))$ points in $\C$ with the following property:
\begin{itemize}
\item
for every $p\in\Pd$, written as $p(z)=c\prod_j(z-\alpha_j)$, there are $d$ points $z^{(1)}, \dots, z^{(d)}\in\Sd$ with $N_p^{\circ n}(z^{(j)})\to \alpha_j$ as $n\to\infty$.
\end{itemize}
Given $\eps>0$, let $n_j\in\N$ be so that $|N_p^{\circ n}(z^{(j)})-\alpha_j|<\eps$ for all $n\ge n_j$. Then the required number of iterations is bounded as follows.
\begin{itemize}  
\item
\textbf{Worst-case complexity:}
we always have 
\[
\sum_j n_j \in 
O\left(d^4\log^2 d+d^3\,\log^2d|\log\eps|\right)
\;.
\] 
\item
\textbf{Average complexity:} if the roots $\alpha_j$ of $p$ are all simple and have mutual distance at least $d^{-k}$ for some $k\in\N$, then 
\begin{equation}
\sum_j n_j\in 
O\left(d^3(\log^2 d)(\log d+k)+d\log|\log\eps|\right) 
\;.
\label{Eq:AverageEstimate}
\end{equation}
In particular, if the roots are randomly distributed in $\disk$, or the coefficients are chosen randomly (subject to the condition that the roots are in $\disk$),
then the expected number of iterations is $O\left(d^3\log^3 d+d\log|\log\eps|\right) $.
\end{itemize}
\end{theorem}

In \cite{NewtonTodor}, we refine these results so that \emph{for randomly distributed roots in $\disk$, one can expect $\sum_jn_j\in O(d^2\log^4d+d\log|\log\eps|)$} (and similar bounds hold if the coefficients are distributed randomly, subject to the restriction that all roots are in $\disk$).
Our current results builds upon earlier work \cite{NewtonFields} that established convergence in polynomial time, but with a rather sub-optimal exponent. 

\medskip

This result measures the complexity in terms of Newton iterations. Of course, each Newton iteration requires arithmetic complexity $d$ (at least the $d$ coefficients of $p$ have to be processed, unless the polynomial is given in special form), but the evaluation of a given polynomial $p$ of degree $d$ at $d$ different points simultaneously is possible (at least when the polynomial is evaluated in terms of coefficients) with arithmetic complexity $O(d\log^2 d)$ using Fast Fourier Transform methods \cite{MoenckBorodin}, \cite[Section~8.5]{HopcroftUllman}. Therefore, the arithmetic complexity differs from the complexity in terms of Newton iterations only by a factor of $\log^2 d$. 

\medskip

\noindent
\emph{Note on parallelization}.
A parallel computer, or a multi-core computer, can take advantage of the inherent parallel structure of the independent Newton iterations, so the algorithm is almost ideally parallelizable. On the other hand, as just mentioned, if polynomials are given in coefficient form, then a single core computer (in classical von Neumann architecture) can compute $d$ independent orbits almost as fast as a single orbit.

\medskip

%\newpage

\noindent
\emph{Polynomial Root-Finding.}
There is an enormous literature on polynomial root-finding; see for instance they surveys by McNamee and Pan\cite{McNamee,McNameeBook,McNameePanBook} and the references therein. Newton's method has been considered difficult to analyze: for instance, Pan \cite{Pan} writes ``Theoretically, the weak point of these algorithms is their heuristic character. \dots \ 
Moreover, in spite of intensive effort of many researchers, convergence of these algorithms has been proved only in the cases where the initial point is already close to a zero or where another similar condition is satisfied.''
In \cite{GLSY}, it is discussed how to discover how the roots are located in the form of clusters, and the difficulty is expressed as ``Then, in the case of a cluster with positive diameter, when arriving close to the cluster, it is well known that the [Newton] iteration may behave in a chaotic way.'' 

The purpose of our work is to control this ``chaotic'' dynamics and to show that the classical simple and stable (and elegant) Newton method is far more efficient than anticipated.
\medskip

\noindent
\emph{Newton's Method in Practice.}
We performed a number of tests on polynomials of large degrees up to $2^{27}>134\cdot 10^6$ jointly with Robin Stoll \cite{NewtonRobin,NewtonRobin2}. In these first set of tests, our method found all roots completely and easily, requiring between $3d^2$ and $6d^2$ iterations to find all roots. For one of these polynomials, there is a sample implementation of MPSolve 3.0, and here on the same computer our algorithm was significantly faster (by orders of magnitude; in part due to the fact that our method allows us to take advantage of the special form of polynomials). More details on practical experiments can be found in \cite{NewtonRobin}. Additional substantial  improvements implemented afterwards made it possible, based on the theory developed in this manuscript, to find all roots of certain degree $2^{20}$ polynomials in about two minutes on a standard PC, and the complexity seems to scale for particular polynomials like $O(d\log d)$ even for degrees up to many millions \cite{NewtonRobin2}. We demonstrate this in Section~\ref{Sec:Experiments}.

One issue that we do not discuss in this paper is to turn the Newton method into an explicit algorithm, including precise stopping criteria and a declaration on the multiplicities of the roots found. This can be done based on our methods; it is more of a technical, rather than a conceptual issue. For details, see \cite{NewtonAlgorithm}.
\medskip

\noindent
\emph{Notation.}
Throughout this text, we will fix a polynomial $p\in\Pd$ of degree $d\ge 2$, and we write $p(z)=c\prod_j(z-\alpha_j)$ and $N_p(z)=z-p(z)/p'(z)$. The coefficient $c$ cancels for Newton's method and will  be omitted. Each root $\alpha_j$ has its \emph{basin} $\hat U_j\subset\C$: this is the set of points that converge to $\alpha_j$ under iteration of $N_p$. The \emph{immediate basin} $U_j$ is the connected component of $\hat U_j$ containing the root $\alpha_j$. As long as we focus attention on a single root, we call it $\alpha$ and its immediate basin $U$. It is well known that $U$ is simply connected and unbounded \cite{Prz} (see also \cite{Mitsu,NewtonInventiones}). Denote by $d_U$ the distance with respect to the unique hyperbolic metric on $U$ with constant curvature $-1$.  Let $D_r(a):=\{z\in\C\colon |z-a|<r\}$ for $a\in\C$ and $r>0$, and let $\disk:=D_1(0)$ be the complex unit disk. Finally, $\log$ always denotes the natural logarithm; sometimes we use the dyadic logarithm and denote it $\log_2$.

\medskip
\noindent
\emph{Overview of the arguments and structure of the paper.}
The first new ingredient in this paper is the concept of ``$R$-central orbits'': these are orbits under Newton's method that stay in the disk $D_R(0)$, so we can maintain control. We will show how to find starting points of Newton's method that are in immediate basins and that have $R$-central orbits. In order to estimate the possible number of iterations, the fundamental idea is the \emph{area used per iteration step}. 
We partition $D_R(0)$ into domains $S_k$ so that for $z\in S_k$, the nearest root has distance approximately $2^{-k}$ from $z$. This will imply that $|z-N_p(z)|\ge \ell:=2^{-k}/d$. We will have orbits $(z_n)$ in the immediate basin $U$ of $\alpha$ for which the hyperbolic distance $d_U(z_n,z_{n+1})= O(\log d)$. Roughly speaking, Euclidean distance bounded below by $\ell$ and hyperbolic distance bounded above by $\log d$ means that $z_n$ and $z_{n+1}$ can be connected by a hyperbolic geodesic segment $\gamma\subset U$ that has Euclidean distance at least $\ell/\log d$ from the boundary, so this curve ``uses up'' an area of approximately $|A_{n,k}(\ell)|\ge \ell^2/\log d=4^{-k}/d^2\log d$ (length times width of the neighborhood of the curve). But $S_k$ is contained in the union of $d$ disks of radius $2^{-k}$ and with total area at most $\pi d4^{-k}$, so there is room for no more than $(d4^{-k})/(4^{-k}/d^2\log d)=d^3\log d$ iterations in each $S_k$ (always up to bounded factors). In the worst case, when there are multiple or near-multiple roots, we will show that we need to consider $k\le O(d)$, hence a total of $O(d^4\log d)$ iterations is required. If the roots are well separated from each other, for instance if they are randomly distributed, it turns out that $k\le\log d$ will usually suffice until the domain of quadratic convergence is reached where $\log|\log\eps|$ further iterations yield precision $\eps$, so we need $O(d^3\log^2d+\log|\log\eps|)$ iterations. All these count the number of iterations required to find a single root. But since all roots are competing for the area, the number of iterations to find \emph{all} roots satisfies the same bounds (except that the small $\log|\log\eps|$ term acquires a factor $d$).

Of course, all these estimates have to be made precise, and we have to make sure that the domains of area do no overlap, which will introduce additional factors of $\log d$. The paper is structured as follows. In Section~\ref{Sec:R-CentralOrbits} we introduce {$R$-central orbits} and show how to find them. In Section~\ref{Sec:GoodStartingPoints}, we construct an explicit finite set of starting points that contains, for each root, at least one $R$-central orbit $(z_n)$ with $d_U(z_n,z_{n+1})=O(\log d)$ as required. 
In Section~\ref{Sec:AreaPerIteration} we estimate the area needed per single iteration step. In Section~\ref{Sec:AreaOrbit}, we estimate the area needed for each orbit: the main point is to make sure that the pieces of area assigned to each orbit point are disjoint; this will introduce another factor $\log d$. 

It remains to discuss when to stop the Newton iteration. If the roots are well-separated from each other, then the iteration reaches the domain of quadratic convergence, and the necessary stopping criterion will be given in Section~\ref{Sec:StoppingCriterion}. We then bring all arguments together and describe the required number of iterations for ``good'' starting points. The worst case of roots that are not well separated, or possibly even multiple, is treated in Section~\ref{Sec:WorstCase}.  

In a brief final Section~\ref{Sec:Experiments}, we report some numerical experiments, obtained jointly with Robin Stoll, that show that for certain polynomials of degrees many millions all roots can be found by Newton's method with guaranteed success and in a matter of only a few hours on standard PC computers using standard arithmetic. This section supports our claim that Newton's method stands out as a root finding method that combines good theory with remarkable practical usefulness.

In Appendix~\ref{Sec:GeometryHypGeodesics}, we prove a general lemma on the area of certain neighborhoods of hyperbolic geodesics in Riemann domains; a major concern is to make sure that these neighborhoods do not overlap. 

Note that in general we are not interested in optimizing constant factors. At a number of places, we specify explicit constants when it is easy to do so, rather than referring to unspecified values. It is certainly possible to improve most constants significantly, but we want to make clear that all of them have very feasible values.

\subsection*{Acknowledgements} I would like to thank especially Todor Bilarev, John Hubbard, Malte Lackmann, Robin Stoll, and Michael Thon for many helpful discussions; I am also grateful to Magnus Aspenberg, Walter Bergweiler, Edward Crane, Bj\"orn Gustafsson, Curt McMullen, Mary Rees, Marcel Oliver, Steffen Roh\-de, Michael Stoll, and an anonymous referee. This research was partially supported by the EU Research and Training Network CODY, by the ESF network HCAA, and by the German Research Council DFG.

%\newpage

\section{Channels and $R$-Central Orbits}
\label{Sec:R-CentralOrbits}

After a brief review of the geometry of immediate basins outside of $\disk$, our main goal in this section is to give a condition on orbits that always stay within a certain Euclidean disk around the origin; such orbits will be called ``$R$-central''.

\begin{figure}
\framebox{\includegraphics[scale=0.8]{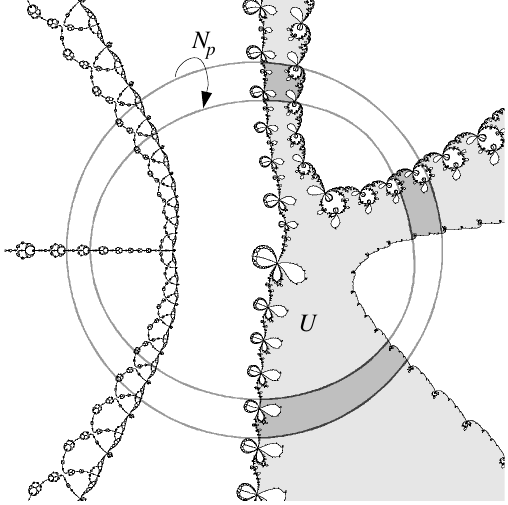}}
\caption{The dynamical plane of a Newton map (of a polynomial of degree $6$). The basin boundaries are black, one immediate basin $U$ with three channels is shaded, and a fundamental annulus for the Newton dynamics is shown. (Picture from \cite{NewtonInventiones}.)}
\label{Fig:GlobalDynamics}
\end{figure}

Consider the immediate basin $U$ of a root $\alpha$ (see Figure~\ref{Fig:GlobalDynamics}); as mentioned above, it is simply connected. The geometry of these basins outside of $\disk$ has been studied in \cite[Section~3]{NewtonInventiones}; in particular, $U$ is unbounded. A \emph{channel} of $U$ is an unbounded connected component of $U\sm\diskbar$, and an \emph{access to $\infty$} of $U$ is a homotopy class (with endpoints fixed) of curves in $U\cup\{\infty\}$ connecting $\alpha$ to $\infty$. Every access to $\infty$ is fixed by $N_p$, and so is every channel: if $B$ is a channel of $U$, then $N_p(B)\sm\diskbar=B$. Each channel contains one access to $\infty$, and each access to $\infty$ defines one channel through which it runs to $\infty$.

Every channel has an associated modulus: the quotient of $B$ by the dynamics of $N_p$ is a conformal annulus $A:=B/(N_p)$ with some modulus $\mu=\mu(A)=:\mu(B)$. One of the main results in \cite{NewtonInventiones} is that each root of a polynomial of degree $d$ has a finite positive number of channels, and at least one of them has modulus $\mu\ge\pi/\log d$ \cite[Propositions~6 and 7]{NewtonInventiones}. 

On the set $B$, we will use three different hyperbolic metrics (all with constant curvature $-1$): $d_B$ is the hyperbolic metric on $B$, while $d_U$ is the hyperbolic metric on $U$ restricted to $B$. Finally, the quotient annulus $A=B/(N_p)$ has a hyperbolic metric, and the infinitesimal metric on $A$ lifts to an infinitesimal metric on $B$, called $d_A$. (Only for $d_B$ is $B$ a complete metric space; for $d_A$ and $d_U$, certain boundary points can be reached in finite distance.) Asymptotically near $\infty$, all three metrics coincide: it is a general principle that if $U_0$ is a hyperbolic Riemann surface and $U_2\subset U_1\subset U_0$ are open subsets, then \emph{restricted to $\ovl{U_2}$ }the hyperbolic metrics of $U_0$ and $U_1$ differ little if $U_0\sm U_1$ is far away from $\ovl{U_2}$ (with respect to the hyperbolic metric of $U_0$); see for instance \cite[Proposition~3.4]{HelenaLasse}. This implies that on all of $B$ we have $d_B>d_A> d_U$: the map $N_p$ is an isometry on $A$ and a contraction on $U$, hence $d_A>d_U$ (for a point $z\in B$, there is a $z_n\in B$ near $\infty$ with $N_p^{\circ n}(z_n)=z$, and asymptotic equality of metrics near $z_n$ implies the claimed inequality). Similarly, there exists a branch $N_p^{-1}\colon B\to B$; this branch contracts $d_B$ (inclusion of hyperbolic domains) and preserves $d_A$, hence $d_B>d_A$ (for $z\in B$, there is a point $z_n:=(N_p^{-1})^{\circ n}$ near $\infty$ for this given branch of $N_p^{-1}$, and the claim follows similarly as above).

For all $z\in B$ with $N_p(z)\in B$, we have $d_A(z,N_p(z))\ge \pi/\mu$ (the core curve in $A$ has hyperbolic length $\pi/\mu$, and every simply closed curve in $A$ is at least as long as the core curve).

Our goal in this section is to prove an existence criterion for orbits that, once they reach $\disk\cap U$, will never leave $D_R(0)$ for a certain radius $R>0$. It is easy to check that if $|z|\ge 1$, then $|N_p(z)|<|z|$ (this follows from \cite[Lemma~3]{NewtonInventiones} as cited above), 
so all orbits outside of $\disk$ move towards $\disk$, and control can be lost only for orbit points in $\disk$ (especially near poles of $N_p$).

%\newpage

\begin{definition}[$R$-central orbits]
\lineclear
An orbit $(z_n)$ will be called \emph{$R$-central} if $|z_n|\le 1$ implies $|z_{n'}|\le R$ for all $n'\ge n$.
\end{definition}

In view of the remark just above, for an orbit to be $R$-central it is sufficient to show that $|z_n|\le 1$ implies $|z_{n+1}|\le R$. 

To locate $R$-central orbits, we need the following definition. 

\begin{definition}[Central subannulus and central channel]
\label{Def:CentralSubannulus} \lineclear
Let $B$ be a channel of an immediate basin $U$, let $A:=B/(N_p)$ be the quotient annulus, and let $\mu:=\bmod A$ be its modulus (so that the core curve of $A$ has length $\pi/\mu$). We define the \emph{central subannulus} of $A$ as the set of points $z\in A$ with injectivity radius less than $2\pi/\mu$. Let the \emph{central subchannel} of $B$ be the preimage in $B$ (under the quotient map) of the central subannulus: this is the set of points $z\in B$ with  $d_A(z,N_p(z))< 2\pi/\mu$. 
\end{definition}

\begin{lemma}[Central subchannel]
\label{Lem:CentralSubchannel} \lineclear
If $A$ has modulus $\mu$, then the central subannulus of $A$ is a parallel subannulus of $A$ with modulus greater than $2\mu/3$.
\end{lemma}
\begin{proof}
Setting $h=\mu/2$, the quotient annulus $A$ is conformally equivalent to the horizontal strip $\{z\in\C\colon |\Im\,z|<h\}$ modulo $z\sim z+1$. The infinitesimal hyperbolic metric on the strip is given by
\[
ds=\frac{\pi|dz|}{2h\cos(\pi|y|/2h)} 
\;,
\]
so the length of the simple closed geodesic in $A$ is $\pi/\mu$, and the parallel subannulus is the set of points $x+iy$ for which $|y|$ satisfies a certain upper bound. In particular, if $\cos(\pi|y|/2h)>1/2$, then the injectivity radius is less than $2\pi/\mu$ (the horizontal curves of Euclidean length $1$ at imaginary part $y$ have hyperbolic lengths less than $2\pi/\mu$, and the corresponding geodesics are shorter than this, but longer than $\pi/\mu$). The condition $\cos(\pi|y|/2h)>1/2$ is satisfied if $|y|<2h/3$, so the central subannulus has modulus greater than $2\mu/3$.
\end{proof}

We can now state the main result of this section.

\begin{theorem}[Large central subchannels have $R$-central orbits]
\label{Thm:SubchannelCentral} \lineclear
Let again $U$ be the immediate basin of a root $\alpha$ and let $B\subset U$ be a channel of $U$ with largest modulus. Then all points in the central subchannel of $B$ have $R$-central orbits, for a value of $R$ that will be specified in Proposition~\ref{Prop:R-CentralOrbits}.
\end{theorem}

The idea of the proof is simple: if $z_0$ is in the central subchannel of $B$, then for its orbit $(z_n)$ we have the estimate  $d_U(z_n,N_p(z_n))\le d_U(z_0,N_p(z_0))<d_A(z_0,N_p(z_0))<2\pi/\mu$. If $|z_n|\le 1$, then the bound $d_U(z_n,N_p(z_n))$ implies an upper bound on $|N_p(z_n)|$ and thus makes the orbit $R$-central. We believe that $R=2$ works for all but a few low values of $\deg p$, but to prove this would require more control on the possible shapes of the channels than we can currently provide. The estimate that follow in this section are relatively weak because they have to account for all possible shapes of channels.  

We start by stating a simple and well-known worst-case estimate of hyperbolic distances in a hyperbolic domain.

\begin{lemma}[Standard bound on hyperbolic arc length]
\label{Lem:BoundArcLength} \lineclear
Let $V$ be a Riemann domain and $p,q\in V$. If $a\in\partial V$ with $s=|p-a|$, then 
\begin{align}
d_V(p,q)\ge \frac 1 2 \int_0^{|p-q|} \frac{dt}{s+t} = \frac 1 2 \log\left(1+\frac{|p-q|}{s}\right)
\;.
\label{Eq:BoundArcLength}
\end{align}
\end{lemma}
\begin{proof}
Infinitesimal hyperbolic distance in $V$ is at least half infinitesimal Euclidean distance, divided by Euclidean distance to $\partial V$, denoted $\operatorname{dist}(\cdot,\partial V)$. If $\gamma\colon[0,|p-q|]\to V$ is a smooth curve parametrized by Euclidean arc length, then its hyperbolic length in $V$ is at least
\[
\frac 1 2 \int_0^{|p-q|} \frac{dt}{\operatorname{dist}(\gamma(t),\partial V)} \ge \frac 1 2 \int_0^{|p-q|} \frac{dt}{\operatorname{dist}(\gamma(0),\partial V)+t}
\;.
\]
In particular, if $\gamma\colon [0,T]\to V$ is the hyperbolic geodesic connecting $p$ to $q$, parametrized by Euclidean arc length $T\ge|p-q|$, then its restriction to $[0,|p-q|]$ has hyperbolic length as least as in \eqref{Eq:BoundArcLength}.
\end{proof}

%\newpage

\begin{lemma}[Hyperbolic distance across fundamental domain]
\label{Lem:HypDistAnnulus} \lineclear
Let $U$ be an immediate basin and $B\subset U$ be a channel of modulus $\mu$. Then all $w,\tilde w\in B$ with $|w|(d-1)/d>|\tilde w|>5$ satisfy $d_U(w,\tilde w)>2/5(\mu+\pi)$.
\end{lemma}
This result is rather weaker than expected: one would expect approximately $d_U(w,\tilde w)\ge \pi/\mu $ (and perhaps a simpler proof), but channels may have complicated geometry; and while some intermediate estimate become less elegant, the final result will be affected only marginally.

\begin{proof}
By \cite[Lemma~3]{NewtonInventiones}, we have $|N_p(w)-w(d-1)/d)|<1/d$ whenever $|w|>1$. 
Let $w':=N_p(w)$; we have $d_B(w,w')>d_A(w,w')\ge \pi/\mu$ and $|w-w'|<(|w|+1)/d$. If all points on $[w,w']$ had Euclidean distance to $\partial V$ at least $2\mu(|w|+1)/\pi d$, then the standard estimate on hyperbolic distance would imply $d_B(w,w')\le 2|w-w'|/(2\mu(|w|+1)/\pi d)<\pi/\mu$, a contradiction. Therefore, there are points $a\in\partial B$ and $w''\in[w,w']$ (the straight line segment from $w$ to $w'$) with $|a-w''|< 2\mu(|w|+1)/\pi d$. This implies
\[
|a-w|\le |a-w''|+|w-w'| <(2\mu/\pi+1)(|w|+1)/d
\;.
\]

The point $a\in\partial B$ is either on $\partial U$ or on $\partial B\cap U$; the latter implies $|a|=1$. 
Suppose first that $a\in\partial U$. Then by Lemma~\ref{Lem:BoundArcLength}
\begin{align*}
d_U(w,\tilde w) 
&\ge\frac 1 2 \log\left(1+\frac{|w-\tilde w|}{(2\mu/\pi+1)(|w|+1)/d}  \right)
\\
&> \frac 1 2 \log\left( 1+\frac{|w|/d}{(2\mu/\pi+1)(|w|+1)/d} \right)
>\frac 1 2 \log\left(1+ \frac{5}{6(2\mu/\pi+1)}\right)
\\
&> \frac{\frac{5}{6(2\mu/\pi+1)}}{2(1+\frac{5}{6(2\mu/\pi+1)})}
%= \frac 1 2 \cdot \frac{5}{6(2\mu/\pi+1)+5} 
= \frac{5\pi}{24\mu+22\pi} > 
 \frac{2}{5(\mu+\pi)}
\;.
\end{align*}

Now we discuss the case that $a\in\partial B\cap U$, i.e.\ $|a|=1$ and 
$|w|-1\le |a-w|< (2\mu/\pi+1)(|w|+1)/d$, hence 
\begin{equation}
\mu>\frac{\pi}{2}\left(d\frac{|w|-1}{|w|+1}-1\right)
> \frac{\pi}2 \left(\frac{4d-6}{6}\right)=\frac{\pi(2d-3)}{6}
\label{Eq:Bound_mu_vs_d}
\end{equation}
or $d<3(\mu/\pi)+3/2$. 
In this case, there is certainly a point $a'\in\partial U$ with $|w-a'|<| w|+1$ (some point $a'\in\disk$) and we get
\begin{align*}
d_U(w,\tilde w)
\;& \ge\; %\int_{0}^{|w-\tilde w|} \frac{dt}{2(| w|+1+t)}=
\frac 1 2 \log\left(1+\frac{|w-\tilde w|}{| w|+1}\right)
\\
&>\;
\frac 1 2 \log\left(1+\frac{| w|/d}{| w|+1}\right)
\ge \frac 1 2 \log\left(1+5/6d\right) > \frac{5/6d}{2(1+5/6d)}
\\
&= \; \frac{5}{12d+10} > \frac{5}{36\mu/\pi+28}
= \frac{5\pi}{36\mu+28\pi} > \frac{2}{5(\mu+\pi)}
\;,
\end{align*}
so the claimed inequality holds in both cases.
\end{proof}

\begin{proposition}[Existence criterion for $R$-central orbits]
\label{Prop:R-CentralOrbits} \lineclear
Suppose all channels of an immediate basin $U$ have modulus at most $\mu$ with $\mu\ge \pi/\log d$, and $z\in U$ has $d_U(z,N_p(z))\le 2\pi/\mu$. Then the orbit of $z$ is $R$-central for $R\le 5(d/(d-1))^{\lceil 5\pi(\log d+1)\rceil}$, and the same holds for all points on the hyperbolic geodesic in $U$ connecting $z$ to $N_p(z)$. 
\end{proposition}
\begin{proof}
Any two points $w,\tilde w\in U$ with $|w|(d-1)/d>|\tilde w|>5$
have $d_U(w,\tilde w)>2/5(\mu+\pi)$ by Lemma~\ref{Lem:HypDistAnnulus}, so going in $U$ from radius $r\ge 5$ to radius $r(d/(d-1))^{k}$ (for some $k\ge 1$) one needs to traverse at least $k$ disjoint concentric annuli with boundary radii differing by a factor of $d/(d-1)$ and hence one needs to traverse at least $k$ complete fundamental domains, and the hyperbolic distance is at least $2k/5(\mu+\pi)$. In particular for $k=\lceil 5\pi(\log d+1)\rceil$ the required hyperbolic distance is at least
\[
\big\lceil5\pi(\log d+1)\big\rceil\cdot \frac{2}{5(\mu+\pi)
} 
\ge 2\pi \frac{\log d+1}{\mu+\pi} 
>
\frac{2\pi}{\mu}
> d_U(z,N_p(z)) \;,
\]
which exceeds the available hyperbolic distance along the orbit of $z$, so this orbit is indeed $R$-central.

Now consider any $z'$ on the hyperbolic geodesic in $U$ connecting $z$ to $N_p(z)$. Then 
\begin{align*}
d_U(z',N_p(z'))
&\le d_U(z',N_p(z))+d_U(N_p(z),N_p(z')) \\ &\le d_U(z',N_p(z))+d_U(z,z')=d_U(z,N_p(z))<2\pi/\mu
\end{align*} 
and the arguments given above also apply to $z'$.
\end{proof}

\begin{remark}
This result provides an upper bound for $R$ that is uniform in $d$: we have $R\le 5(d/(d-1))^{\lceil5\pi(\log d+1)\rceil} < 5e^{\lceil5\pi(\log d+1)\rceil/d}$. More precisely, we have $R<50$ for $d>30$, and $R<10$ for $d>133$. Ultimately, the precise value of $R$ is not of too large importance as it will enter our estimates only logarithmically: the number of iterations scales with $M\in\N$ such that $2^M-1\ge R$ (see Proposition~\ref{Prop:IsolatedRoots}); we can use $M=5$ for $d\ge 41$ and $M=3$ for $d\ge 316$. Presumably, $R\le 3$ and $M=2$ work for all but very few values small of $d$.
\end{remark}

\begin{proof}[Proof of Theorem~\ref{Thm:SubchannelCentral}]
Let $\mu$ be the modulus of $B$. If $z$ is in the central subchannel of $B$, then $d_U(z,N_p(z))<d_A(z,N_p(z))<2\pi/\mu$. By Proposition~\ref{Prop:R-CentralOrbits}, the orbit of $z$ is $R$-central.
\end{proof}

\section{Good Starting Points}
\label{Sec:GoodStartingPoints}

In this section, we construct a finite set of starting points $\Sd$ depending only on the degree $d$ (and the normalization of $p$) so that for each root $\alpha$ of $p$, one of the points $z\in\Sd$ is in the immediate basin $U_\alpha$ and satisfies $d_U(z,N_p(z))\le 2\log d$, and so that the orbit of $z$ is $R$-central for $R$ as in Proposition~\ref{Prop:R-CentralOrbits}.

All we need to do is specify a finite set of starting points that will intersect, for every root, the central subchannel of the channels with largest modulus. Since every root has a channel with modulus at least $\pi/\log d$, it is sufficient to specify a finite set of starting points that intersects all subchannels of all channels with moduli at least $2\pi/3\log d$. This can be accomplished by the methods in \cite{NewtonInventiones}, so we can now construct an explicit point grid (in that paper, we used $.2663\log d$ concentric circles that each contain $8.33d\log d$ points; here we use $(3/2)$ as many circles because we want to hit the channel within the central subchannel with $(2/3)$ the modulus).

%\newpage

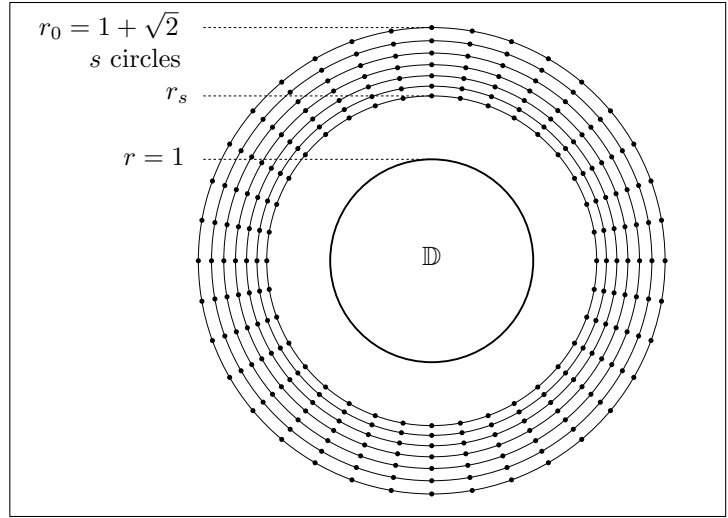
\begin{figure}[htbp]
\begin{picture}(0,200)(130,0)
\framebox{\rule{17mm}{0pt}\includegraphics[width=.6\textwidth,trim=0 30 0 28]{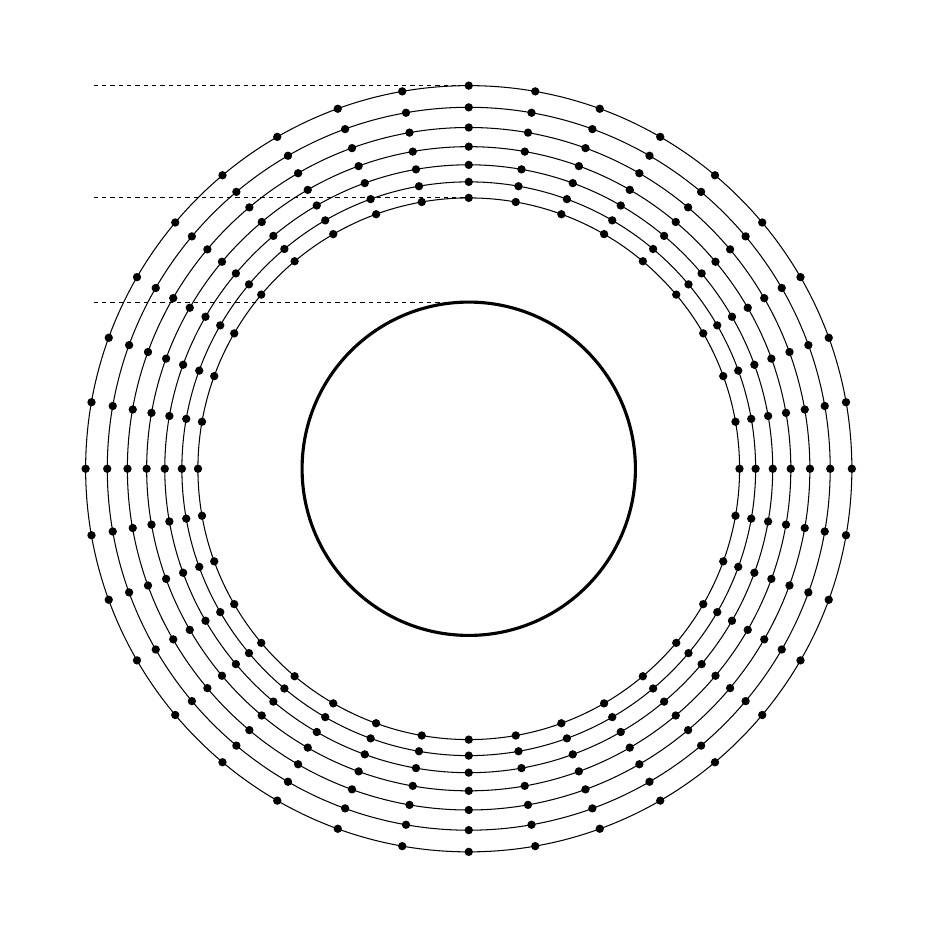}}
\put(-228,130){$r=1$}
\put(-212,155){$r_s$}
\put(-241,167){$s$ circles}
\put(-260,179){$r_0=1+\sqrt 2$}
\put(-115,92){$\disk$}
\end{picture}
\caption{Schematic sketch of the required starting points: they are equidistributed on $s=\lceil 0.4\log d\rceil$ concentric circles around the unit disk that contains all roots. }
\label{Fig:StartingPointSketch}
\end{figure}

\begin{definition}[Efficient grid of starting points]
\label{Def:EfficientGrid} \lineclear
For each degree $d$, we construct a circular grid $\Sd$ of starting points as follows (as sketched in Figure~\ref{Fig:StartingPointSketch}).
For $\nu=1,2,\dots,s=\lceil 0.4 \log d\rceil$, set 
\[
r_\nu:=(1+\sqrt 2)\left(\frac{d-1}{d}\right)^{(\nu-1/2)/2s}
\] 
and for each circle around $0$ of radius $r_\nu$, choose $\lceil 8.33 d\log d\rceil$ equidistant points (independently for all the circles).
\end{definition}

\begin{theorem}[Efficient grid of starting points]
\label{Thm:EfficientGrid}\lineclear
For each degree $d$, the set $\Sd$ contains $3.33 d\log^2 d (1+o(1))$ points, and it has the following universal property. If $p$ is any complex polynomial of degree $d$, normalized so that all its roots are in $\disk$, then there are $d$ points in $\Sd$ that converge to the $d$ roots of $p$, so that 
for each root, one of these points is in the central subchannel of a channel with largest modulus. The corresponding orbits are all $R$-central for a uniform value of $R$ (as in Proposition~\ref{Prop:R-CentralOrbits}). More precisely, if $\mu$ is the largest modulus of all channels of a root $\alpha$ with immediate basin $U$, then the corresponding orbit converging to $\alpha$ satisfies  $d_U(z,N_p(z))<2\pi/\mu<2\log d$.
\end{theorem}
\begin{proof}
The annulus 
\[
V:=\left\{z\in\C\colon  (1+\sqrt 2)\sqrt{(d-1)/d}<|z|<(1+\sqrt 2)\right\}
\] 
is contained in a fundamental domain of the Newton dynamics by \cite[Lemma~10]{NewtonInventiones}. For $\nu=1,2,\dots,s=\lceil 0.4\log d\rceil$, subdivide $B$ into $s$ subannuli
\[
V_\nu=\left\{z\in\C\colon (1+\sqrt 2)\left(\frac{d-1}{d}\right)^{\nu/2s} < |z| < (1+\sqrt 2)\left(\frac{d-1}{d}\right)^{(\nu-1)/2s}
\right\}
\;.
\]

Since $V$ is contained in a fundamental domain of the dynamics, each subchannel of any root with modulus $\mu$ intersects $V$ in a quadrilateral of modulus at least $\mu$, and by the Gr\"otzsch inequality, it intersects at least one $V_\nu$ in a quadrilateral with modulus at least $s\mu$.
Each root has a channel with modulus at least $\pi/\log d$, so the central subchannel has modulus at least $2\pi/3\log d$, and this central subchannel intersects some $V_\nu$ in a quadrilateral with modulus at least $2s\pi/3\log d > 0.2663\pi$, independent of $d$. In \cite[Section~6]{NewtonInventiones}, it is shown that $\lceil 8.3254 d\log d\rceil$ equidistributed points on each of these circles will find all quadrilaterals connecting the boundaries of the $B_\nu$ with modulus at least $0.2663\pi$, as in our case.

Therefore, the grid $\Sd$ intersects all central subchannels of all largest channels, and the claim follows.
\end{proof}

\begin{remark}
The number of starting points of $O(d\log^2d)$ from \cite{NewtonInventiones} has been further reduced to $O(d(\log\log d)^2)$ in \cite{BLS}, by using a probabilistic set of starting points. This approach could also be used in our case.
\end{remark}

Now we have a good set of starting points leading to $R$-central orbits. We proceed to estimate the number of required iterations.

%\newpage

\section{Area per Iteration Step}
\label{Sec:AreaPerIteration}

In this section, we show that every iteration step ``uses up'' a certain area in the plane; since $R$-central orbits remain within some disk $D_R(0)$, this will provide an upper bound on the possible number of iterations.

Consider some point $z\in U$, set $\tau:=d_U(z,N_p(z))$, and let $\gamma\colon[0,T]\to U$ be the hyperbolic geodesic connecting $z$ to $N_p(z)$, parametrized by \emph{Euclidean} arc length. For $t\in[0,T]$, let $\eta(t)$ be the Euclidean distance from $\gamma(t)$ to $\partial U$, and let $X(t)$ be the
Euclidean straight line segment of length $\eta(t)$ (not containing the endpoints) with center at $\gamma(t)$ and perpendicular to $\gamma$ at $\gamma(t)$, so that $\gamma(t)$ disconnects $X(t)$ into two open segments of length $\eta(t)/2$; see Figure~\ref{Fig:X(t)}. All the segments $X(t)$ are disjoint (Proposition~\ref{Prop:All_Xt_Disjoint} in the appendix). 

For $k\in\N$, let $X_{k}(t)$ be the restriction of $X(t)$ to lengths at most $2^{-k+1}$ (i.e., $X_{k}$ is the perpendicular line segment to $\gamma(t)$ centered at $\gamma(t)$ and extending in both directions for a length of $\min(\eta(t)/2,2^{-k})$\,). 
Let $A_{k}:=\bigcup_{t\in[0,T]} X_k(t)$ be the subset of $U$ covered by the $X_{k}(t)$ for $t\in[0,T]$, and for $\ell\le T$ let $A_k(\ell):=\bigcup_{t\in[0,\ell]} X_k(t)$ be the analogous set for $t\in[0,\ell]$. This set of course depends on $z$, so a more explicit description of this set would be $A_k(\ell,z)$. We denote the Euclidean area of $A_k(\ell)$ by $|A_k(\ell)|$.
Of course, $A_k(\ell)\subset A_{k-1}(\ell)$, and the limit as $k\to-\infty$ is $A_{-\infty}(\ell)=\bigcup_{t\in[0,\ell]}X(t)$. 

\begin{figure}[htbp]
\framebox{
\includegraphics[scale=.7]{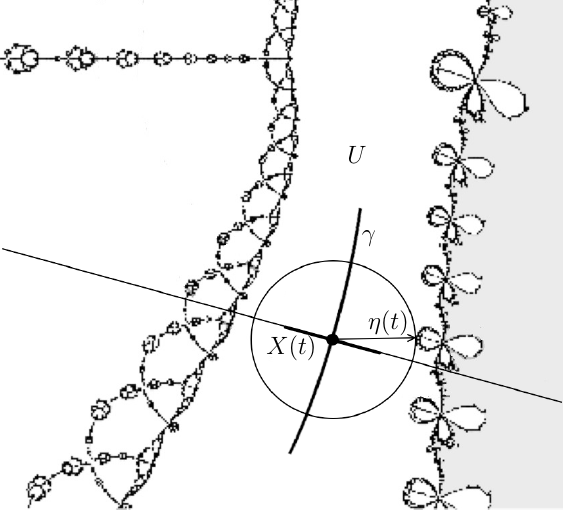}
}\caption{A line segment $X(t)$ perpendicular to the geodesic $\gamma$; in both directions from the point $\gamma(t)$, it extends to a length of $\eta(t)/2$, where $\eta(t)$ is the distance from $\gamma(t)$ to $\partial U$. (Based on a picture from \cite{NewtonInventiones}.)}
\label{Fig:X(t)}
\end{figure}

\begin{lemma}[Area per iteration step]
\label{Lem:AreaPerIteration} \lineclear
We have $|A_{k}(\ell)|\ge {\ell^2}/{(2\tau+2^{k-1}\ell)} $.
Similarly, the set $\bigcup_{t\in[0,\ell]} X(t)$ has Euclidean area at least $\ell^2/2\tau$.
\end{lemma}
\begin{proof}
The hyperbolic length of $\gamma$ between $z$ and $N_p(z)$ satisfies
\[ 
\tau\ge \frac 1 2 \int_0^{T} \frac{dt}{\eta(t)}
\ge
 \frac 1 2 \int_0^{\ell} \frac{dt}{\eta(t)}
\;.
\]
This implies
\[
\int_0^{\ell} \frac{dt}{\min\left(\eta(t),2^{-{k+1}}\right)}
\le 
\int_0^{\ell} \frac{dt}{\eta(t)} + \int_0^{\ell} \frac{dt}{2^{-k+1}}
\le 2\tau + 2^{k-1} \ell \;.
\]

Since all $X_k(t)$ are disjoint, an exercise in elementary differential geometry shows that
\[
|A_{k}(\ell)|= \int_0^{\ell} 2\min\left(\eta(t)/2,2^{-k}\right)\,dt
= \int_0^{\ell} \min\left(\eta(t),2^{-k+1}\right)\,dt
\;.
\]

The Cauchy-Schwarz inequality for the functions $\sqrt{\min\left(\eta(t),2^{-k+1}\right)}$ and $1/\sqrt{\min\left(\eta(t),2^{-k+1}\right)}$ implies
\begin{align*}
|A_{k}(\ell)| &= \int_0^{\ell} \min\left(\eta(t),2^{-k+1}\right) \,dt \ge
\frac{\ell^2}{\displaystyle\int_0^{\ell} \frac{1}{\min\left(\eta(t),2^{-k+1}\right)}\,dt}
\\
&\ge
\frac{\ell^2}{2\tau+2^{k-1}\ell} 
\;. 
\end{align*}
The computation for $\bigcup_{t\in[0,\ell]} X(t)$ is even simpler, replacing the $\min$ by $\eta(t)$.
\end{proof}

Define the following subsets of $\C$, for $k\in\Z$:
\[
S_k^*:=\left\{z\in\C\colon \min_{j} |z-\alpha_j| \le 2^{-k} \right\} 
\]
and
\[
S_k:=\left\{z\in\C\colon \min_{j} |z-\alpha_j| \in \left(2^{-(k+1)},2^{-k}\right] \right\}
= S_k^* \sm S_{k+1}^* \,\,.
\]
Every $S^*_k$ the union of closed disks of radius $2^{-k}$ around all the roots, and for every $M\in\Z$ the disk $D_{2^M-1}(0)$ is partitioned by the $S_k$ for $k\ge -M $ (appropriately restricted).

\begin{lemma}[Distance from $A_k(\ell)$ to roots]
\label{Lem:AreaInDisk} \lineclear
If $z\in S^*_k$ and $|z-N_p(z)|\ge\ell$, then each point in $A_k(\ell)$ has Euclidean distance at most $2^{-k+1}+\ell$ from a root of $p$.
\end{lemma}
\begin{proof}
Since $z\in S^*_k$, the point $z$ has distance at most $2^{-k}$ from a root; then in $A_k(\ell)$ we consider only a segment of length $\ell$ on the geodesic from $z$ to $N_p(z)$, and each point in $A_k(\ell)$ has distance at most $2^{-k}$ from a point on this geodesic segment. 
\end{proof}

%\newpage

\section{Area Along an Orbit}
\label{Sec:AreaOrbit} 

%\newpage

In this section, we investigate how much area is ``used up'' along an orbit $(z_n)$: this involves discussing when the area used for different iteration steps is disjoint, and determining a geometric condition on the hyperbolic displacement in each Newton step before the orbit enters the domain of quadratic convergence near a root; this condition is used for estimating how many iterations the orbit can spend within each $S_k$.

We continue to consider the immediate basin $U$ of some root $\alpha$; in this section, we will assume this root is simple. 
Let $\phi\colon U\to\disk$ be a Riemann map with $\phi(\alpha)=0$; then $f:=\phi\circ N_p\circ\phi^{-1}\colon\disk\to\disk$ is holomorphic with $f(0)=0=f'(0)$. Hence $f(w)/w^2$ is still holomorphic and sends $\disk$ to itself, and it follows that $|\phi(N_p(z))|=|f(w)|\le |w|^2=|\phi(z)|^2<|\phi(z)|$. 

Smale \cite{Smale} has the concept of ``approximate zero'': that is a point sufficiently near a simple root from which the convergence is quadratic. We begin with a lemma that gives useful dynamic consequences both for the case when a finite orbit point is (not yet) an approximate root (part (1)) and when it (part (2)): for us, a practical criterion is whether or not $|\phi(z)|<e^{1/2}-1$.

Consider the orbit $z_n=N_p^{\circ n}(z_0)$ for a point $z_0\in U$ and let $\tau:=d_U(z_0,z_1)$; then $d_U(z_n,z_{n+1})\le \tau$ for all $n$. (Eventually, we will use $\tau=O(\log d)$.)

\begin{lemma}[Hyperbolic distance to root and Newton dynamics]
\label{Lem:HypDistanceOrbit} \lineclear
(1) If $n'>n$ so that $|\phi(z_n)|> |\phi(z_{n'})|\ge e^{1/2}-1\approx 0.649$, then 
\[
d_U(z_n,z_{n'})\ge d_\disk(|\phi(z_n)|,|\phi(z_{n'}|)> (n'-n)/2
\;;
\]

\noindent
(2) if $|\phi(z_{n})|<1/2$, then $|z_{n'}-\alpha|<\eps$ for all $n'> n+\log_2 |\log_2\eps -5|$, and 

\noindent
(3) if $|\phi(z_{n})|<e^{1/2}-1$, then $|z_{n'}-\alpha|<\eps$ for all $n'> n+1+\log_2 |\log_2\eps -5|$.

\end{lemma}
\begin{proof}
If $|w|\ge e^{1/2}-1$, then 
\begin{eqnarray*}
d_\disk(|w|,|f(w)|)& > &\left|\log(1-|w|)-\log(1-|f(w)|)\rule{0pt}{10.5pt}\right|
\\
&\ge& \left| \log(1-|w|) -\log(1-|w|^2)\right|
= \log(1+|w|) \ge 1/2
\end{eqnarray*}
(in the first inequality, we used the fact that the hyperbolic distance in $\disk$  exceeds the hyperbolic distance in the smallest left half plane containing $\disk$). Therefore, if $|\phi(z)|\ge e^{1/2}-1$, then
\begin{align*}
d_U(z,N_p(z))&=d_\disk(\phi(z),\phi(N_p(z))
\\
&=d_\disk(\phi(z),f(\phi(z))\ge d_\disk(|\phi(z)|,|f(\phi(z)|)\ge 1/2
\;.
\end{align*} 
Recall that $|\phi(z_n)|>|\phi(z_{n+1})|>\dots>|\phi(z'_n)|$; as long as all these absolute values are $e^{1/2}-1$ or greater, every subsequent iteration adds hyperbolic distance at least $1/2$:
\begin{align*}
d_U(z_n,z'_n) &= d_\disk(\phi(z_n),\phi(z_{n'})\ge d_\disk(|\phi(z_n)|,|\phi(z_{n'}|) 
\\ 
&=\sum_{k=n}^{n'-1} d_\disk(|\phi(z_n)|,|\phi(z_{n'}|) \ge (n'-n)/2
\end{align*}
and the first claim follows.

\looseness-1
If $|\phi(z)|=|w|<1/2$, then $|\phi(N_p^{\circ m}(z))|=|f^{\circ m}(w)|\le 2^{-2^m}$. By the Koebe $1/4$-theorem, $|(\phi^{-1})|'(0) < 8$ (there are roots other than $\alpha$ in $\disk$, so not all of $\disk$ can be in $U$). By the Koebe distortion theorem, $|N_p^{\circ m}(z)-\alpha| < 32\cdot 2^{-2^{m}}$, and this is less than $\eps$ provided $m> \log_2|\log_2(\eps/32)| = \log_2 |\log_2\eps -5|$. 

Finally, if $|\phi(z)|<e^{1/2}-1$, then $|\phi(N_p(z))| \le (e^{1/2}-1)^2<1/2$.
\end{proof}

%\goodbreak

Define sets $A_{n,k}(\ell):=A_k(\ell,z_n)$: these are the sets $A_k(\ell)$ based at the points $z_n$, as defined at the beginning of Section~\ref{Sec:AreaPerIteration}.

%\newpage

\begin{lemma}[Disjointness of areas]
\label{Lem:DisjointAreas} \lineclear
The sets $A_{n,k}(\ell)$ and $A_{n',k'}(\ell)$ are disjoint if $n'-n> 2\tau+4\log 3$ and $|\phi(z_{n'})|\ge e^{1/2}-1$.
\end{lemma}
\begin{proof}
Since $A_{n,k}(\ell)\subset A_{n,k-1}(\ell)\subset A_{n,-\infty}(\ell)$, it is sufficient to consider only the case $A_{n,-\infty}(\ell)$, i.e.\ with perpendicular segments $X(t)$ of length $\eta(t)$ independent of $k$.

If $n'> n+2\tau+4\log 3$, then $d_U(z_n,z_{n'})\ge d_\disk(|\phi(z_n)|,|\phi(z_{n'}))> \tau+2\log 3$ by Lemma~\ref{Lem:HypDistanceOrbit}. Since $|\phi(z_{n'+1})|<|\phi(z_{n'})|$ and the disk $D_{|\phi(z_{n'})|}(0)\subset\disk$ is geodesically convex \cite{Jorgensen}, the entire geodesic segment from $\phi(z_{n'})$ to $\phi(z_{n'+1})$ has distance at least $\tau+2\log 3$ from $\phi(z_{n})$. Since the geodesic segment connecting $\phi(z_n)$ to $\phi(z_{n+1})$ has hyperbolic length at most $\tau$, all points on this geodesic segment have distance at least $2\log 3$ from all points on the geodesic segment connecting $\phi(z_{n'})$ to $\phi(z_{n'+1})$. The same holds for the analogous geodesic segments in $U$.

Essentially by definition, each point on $X(t)$ has hyperbolic distance from its midpoint $\gamma(t)$ of less than $\log 3$ (if $V\subset\C$ is a Riemann domain and $D_r(a)\subset V$, then all points in $D_{r/2}(a)$ have hyperbolic distance in $V$ of at most $\log 3$). The claim follows.
\end{proof}
%\newpage

%\newpage

\begin{proposition}[Number of points in $S^*_k$]
\label{Prop:NumberOrbitPoints} \lineclear
For every $k$ and every $\ell>0$, the set $S^*_k$ contains at most 
\begin{equation}
\pi d (2^{-k+1}+\ell)^2(2\tau+2^{k-1}\ell)\lceil 2\tau+4\log 3\rceil\ell^{-2}
\label{Eq:NumberOrbitPointsS}
\end{equation}
orbit points $z_n$ with $|z_n-z_{n+1}|\ge \ell$ and $|\phi(z_n)|\ge e^{1/2}-1$.
\end{proposition}
\begin{proof}
If $z_n\in S^*_k$ and $|z_n-z_{n+1}|\ge \ell$, then each point in $A_{n,k}(\ell)$ has Euclidean distance at most $2^{-k+1}+\ell$ from some root by Lemma~\ref{Lem:AreaInDisk}, so $A_{n,k}(\ell)$ is contained in a set of total area at most $\pi d (2^{-k+1}+\ell)^2$. Each $A_{n,k}(\ell)$ has area at least $\ell^2/(2\tau+2^{k-1}\ell)$ by Lemma~\ref{Lem:AreaPerIteration}, and by Lemma~\ref{Lem:DisjointAreas} the sets $A_{n,k}(\ell)$ and $A_{n',k}(\ell)$ are disjoint if $n'-n> 2\tau+4\log 3$ and $|\phi(z_{n})|\ge e^{1/2}-1$. Therefore, there can be at most
\[
\frac{\pi d (2^{-k+1}+\ell)^2\lceil 2\tau+4\log 3\rceil }{\ell^2/(2\tau+2^{k-1}\ell)}
=
\frac{\pi d (2^{-k+1}+\ell)^2(2\tau+2^{k-1}\ell)\lceil 2\tau+4\log 3 \rceil}{\ell^2}
\]
such points, for any choice of $\ell$. 
\end{proof}

%\newpage

\begin{remark} \looseness-1
The sets $A_{n,k}(\ell)$ are contained in $U$, so in the end the various orbits in the different immediate basins $U_{\alpha}$ for different roots $\alpha$ will compete for the area. The last result can thus be sharpened as follows. For a root $\alpha$, let 
\[
U_{\alpha,k}(\ell):=\{z\in U_{\alpha}\colon|z-\alpha_j|<2^{-k+1}+\ell \mbox{\; for some root $\alpha_j$}\} 
\]
(this is the $2^{-k}+\ell$-neighborhood of $S^*_k$ restricted to $U_\alpha$).
Then the set $S^*_k$ contains at most
\begin{equation}
|U_{\alpha,k}(\ell)|\cdot (2\tau+2^{k-1}\ell)\lceil 2\tau+4\log 3\rceil \ell^{-2}
\label{Eq:AreaSingleImmedBasin}
\end{equation}
points on the orbit $(z_n)\subset U_\alpha$ with $|z_n-z_{n+1}|\ge \ell$ and $|\phi(z_n)|\ge e^{1/2}-1$, and of course we have
\begin{equation}
\sum_j|U_{\alpha_j,k}(\ell)|\le \pi d (2^{-k+1}+\ell)^2
\;.
\label{Eq:SumAreaU_alpha}
\end{equation}
\end{remark}

\begin{lemma}[Newton displacement and nearest root]
\label{Lem:LogDerivative} \lineclear
For any $z\in\C$, the nearest root $\alpha$ satisfies $|z-\alpha|\le d|z-N_p(z)|$.
\end{lemma}
\begin{proof}
This is easy and well known:
\begin{equation}
z-N_p(z) = \frac{1}{p'(z)/p(z)} = \frac{1}{\sum_{\alpha_j}\frac{1}{z-\alpha_j}}
\;,
\label{Eq:NewtonDisplacement}
\end{equation}
hence
\[
\frac{1}{|z-N_p(z)|} \leq \sum_{\alpha_j} \frac{1}{|z-\alpha_j|}
\le d\frac{1}{\inf_{\alpha_j}|z-\alpha_j|} \;.
\qedhere
\]
\end{proof}

\begin{corollary}[Number of points in $S_k$]
\label{Cor:NumberIteratesS_k} \lineclear
For any $k$, the set $S_k$ contains at most 
\begin{equation}
\pi (4d+1)^2(2\tau d+1/4)\lceil 2\tau+4\log 3\rceil \in O(d^3\tau^2)
\label{Eq:NumberOrbitPoints}
\end{equation} 
points on any orbit $(z_n)\subset U$ with $|\phi(z_n)|>e^{1/2}-1$. 
\end{corollary}
\begin{proof}
If $z_n\in S_k$, then $|z_n-z_{n+1}|> 1/d2^{k+1}$ by Lemma~\ref{Lem:LogDerivative}, so we use $\ell=1/d2^{k+1}$ in Proposition~\ref{Prop:NumberOrbitPoints} and obtain the estimate
\begin{eqnarray*}
&&\pi d^34^{k+1}\left(2^{-k+1}+\frac{2^{-k-1}}{d}\right)^2\left(2\tau+\frac 1{4d}\right)\lceil 2\tau+4\log 3 \rceil
\\
&=&
\pi (4d+1)^2(2\tau d+1/4)\lceil  2\tau+4\log 3\rceil
\end{eqnarray*}
as claimed.
\end{proof}

%\goodbreak

\begin{remark}
As before (see \eqref{Eq:AreaSingleImmedBasin}), the different roots have to compete for the total area available, and the set $S_k$ can contain at most
\begin{equation}
|U_{\alpha,k}(1/d2^{k+1})|\cdot  (2\tau+1/4d)\lceil 2\tau+4\log 3\rceil d^24^{k+1} 
\label{Eq:NumberIteratesS_k}
\end{equation}
points on the orbit $(z_n)\subset U_\alpha$ with $|\phi(z_n)|>e^{1/2}-1$.
\end{remark}

%\newpage

\section{Uniformly Separated Roots}
\label{Sec:StoppingCriterion}

Now that we have good bounds on how many iterations any of our selected orbits can spend within each $S_k$, we have to discuss the possible values of $k$. Any disk $D_R(0)$ is partitioned by $S_k\cap D_R(0)$ for $k\ge -\log_2(R+1)$, and we gave an upper bound for $R$, hence a lower bound for $k$, in Proposition~\ref{Prop:R-CentralOrbits}. We also need an upper bound for $k$, that is a ``stopping criterion'' when the orbit is sufficiently close to a root.

We will need two kinds of stopping criteria: a worst-case estimate that applies especially when there are multiple or near-multiple roots, and a better estimate in case the roots are reasonably well separated from each other, so the orbit is already an approximate root. We first investigate well-separated roots: we say that the roots are \emph{$\delta$-separated}  if they are all simple and have mutual distance at least $\delta$. If roots are randomly distributed in $\disk$, with high probability they will be $\delta$-separated with $\delta = O(1/d)$ (see the remark at the end of this section). 
Multiple or near-multiple roots will be treated in Section~\ref{Sec:WorstCase}. 

\begin{lemma}[Stopping criterion]
\label{Lem:StoppingCriterion2} \lineclear
(1) If $|z-\alpha|<|z-\alpha'|/2d$ for all roots $\alpha'\neq\alpha$, then the Newton orbit of $z$ converges to $\alpha$. 

(2)
If even $|z-\alpha|<|z-\alpha'|/(4d+3)$ for all $\alpha'\neq \alpha$, then $|N^{\circ n}(z)-\alpha|<\eps$ for all $n> \log_2|\log_2\eps -5|$.
\end{lemma}

\begin{proof}
(1) We may rescale coordinates by an automorphism of $\C$ so that $z=0$ and $\alpha=1$. By hypothesis, we have $|z-\alpha|=1$ and $|z-\alpha_j|> 2d$ for all $\alpha_j\neq\alpha$. As in the proof of Lemma~\ref{Lem:LogDerivative}, this implies $|\sum_{\alpha_j\neq\alpha} 1/(z-\alpha_j)| < (d-1)/2d<1/2$, so $\sum_{\alpha_j} 1/(z-\alpha_j)\in D_{1/2}(-1)$, the open disk of radius $1/2$ around $-1$. Thus $z-N_p(z) = (\sum_{\alpha_j}1/(z-\alpha_j))^{-1} \in D_{2/3}(-4/3)$ (the image of a circle under $z\mapsto 1/z$ is a circle that in this case is real symmetric, and it is easy to compute the points where it intersects the real line). Therefore $N_p(z)\in D_{2/3}(4/3)$ and $|N_p(z)-\alpha|<1=|z-\alpha|$, so by induction the orbit of $z$ converges to $\alpha$.

(2) If $|z-\alpha|<|z-\alpha'|/(4d+3)$, then again we choose coordinates with $\alpha=1$ and $z=0$, so all $|\alpha'|>4d+3$. All $z'\in D_2(1)$ have $|z'-\alpha|<2< |z'-\alpha'|/2d$, so $D_2(1)$ is contained in the (rescaled) immediate basin of $\alpha$ and we have $d_U(z,\alpha)<d_\disk(1/2,0)$, hence $|\phi(z)|<1/2$. The claim thus follows from Lemma~\ref{Lem:HypDistanceOrbit} (2).
\end{proof} 

%\newpage

We would like to point out that the following result does not require that all roots are $\delta$-separated, but only that we have some root $\alpha$ that has distance at least $\delta$ from all other roots (which may be multiple or clustered).

\begin{proposition}[Number of iterations on orbit, $\delta$-separated case]
\label{Prop:IsolatedRoots} \lineclear
Suppose $\alpha$ is a simple root and $|\alpha'-\alpha|>\delta$ for all roots $\alpha'\neq \alpha$. If $(z_n)$ is an $R$-central orbit in $U$ with $|z_0|\le R\le 2^M-1$ and $d_U(z_0,z_1)\le\tau$, then we have $|z_N-\alpha|<\eps$ for all $N$ at least
\begin{align*}
&\pi (4d+1)^2(2\tau d+1/4)\lceil 2\tau+4\log 3\rceil\lceil \log_2(4d+4)/\delta+M+1\rceil 
\\ &\qquad +\log_2|\log_2\eps-5|
 \\
& = (64\pi d^3 \tau^2 (\log_2d+|\log_2\delta|)+\log_2\log\eps) (1+o(1)) 
 \\
& \in O\left(d^3\tau^2(\log d+|\log\delta|)+\log|\log\eps|\right) \;.
\end{align*}

\end{proposition}
\begin{proof}
If $|z-\alpha|<\delta/(4d+4)$, then all roots $\alpha'\neq\alpha$ satisfy $|z-\alpha'|> \delta-\delta/(4d+4)=(4d+3)\delta/(4d+4)$, so the orbit of $z$ satisfies the hypothesis of Lemma~\ref{Lem:StoppingCriterion2} (2) and will be $\eps$-close to $\alpha$ after at most $\log_2|\log_2\eps-5|$ iterations.

Therefore, choose $K\in\N$ with $2^{-K}\le \delta/(4d+4)$, i.e., $K=\lceil \log_2((4d+4)/\delta)\rceil$. We only have to consider the number of iterations that the orbit stays in $S_k$ with $k\le K$. Since the orbit is contained within $D_{2^M-1}$ by hypothesis, we have $k\ge -M$, so we need to consider $k\in\{-M,-M+1,\dots,K\}$. 

By Corollary~\ref{Cor:NumberIteratesS_k}, any orbit $(z_n)$ with $d_U(z_0,z_1)<\tau$ has at most $\pi (4d+1)^2(2\tau d+1/4)\lceil 2\tau+6\rceil$ points within each $S_k$, so the total number of iterations required for the orbit $z_n$ is at most 
\[
\pi (4d+1)^2(2\tau d+1/4)\lceil 2\tau+6\rceil \lceil \log_2((4d+4)/\delta)+M+1\rceil+\log_2|\log_2\eps-5|
\;.
\]
\end{proof}

Note that here and elsewhere when we give asymptotic complexity results in $O$-notation we always have explicit constants, and these are small (so we are not hiding gigantic constants behind this notation). 

\begin{remark}
Again, the $d$ roots have to compete for the available area within $\disk$. If all roots are simple and $\delta$-separated, and there are $d$ orbits, one in each immediate basin, that satisfy the hypotheses of Proposition~\ref{Prop:IsolatedRoots}, then the combined number of iterations required to reach $\eps$-precision for all $d$ roots is at most
\begin{align*}
\pi (4d+1)^2\left(2\tau d+\frac 1 4\right)\lceil 2\tau+4\log 3\rceil\left\lceil \log_2\left(\frac{4d+4}{\delta}\right)+M+1\right\rceil\\
+d\log_2|\log_2\eps-5| \rule{60mm}{0pt}
 \\
 \in O\left(d^3\tau^2(\log d+|\log\delta|+M)+d\log|\log\eps|\right) \;;
\end{align*}
this is almost the same bound as for a single orbit (each area element can be used for only one root), except that the estimate $\log_2|\log_2\eps-5|$ (which takes care of approximate roots and does not involve area) applies for each root separately. 
\end{remark}

%\newpage

\begin{theorem}[Efficient grid of starting points]
\label{Thm:EfficientGridSeparated}\lineclear
For each degree $d$, the set $\Sd$ has the following universal property. If $p$ is a complex polynomial of degree $d$, normalized so that all its roots are in $\disk$, and so that all roots are simple and have mutual distance at least $\delta$, then the universal starting point set $\Sd$ contains $d$ points that converge to the $d$ roots of $p$ and so that the combined number of iterations required to reach $\eps$-precision is at most 
\begin{align*}
\pi (4d+1)^2(4d\log d+1/4)\lceil 4\log d +4\log 3\rceil\times \\ \times \lceil \log_2((4d+4)/\delta)+M+1\rceil+d\log_2|\log_2\eps-5|
 \\
= \left(256\pi d^3 \log^2 d (\log d+|\log\delta|) +d\log_2\log\eps\right)(1+o(1))
\\
\in O\left(d^3\log^2 d(\log d+|\log\delta|)+d\log|\log\eps|\right) 
\end{align*}
where $M$ is such that $2^M-1\ge R$ from Proposition~\ref{Prop:R-CentralOrbits}.
If not all roots are $\delta$-separated, then a subset of these $d$ points finds all those roots that are $\delta$-separated from all other roots with $\eps$-precision in the given number of iterations.
\end{theorem}
Recall from the remark after Proposition~\ref{Prop:R-CentralOrbits} that $R$ satisfies an explicit bound for every $d$, and is universally bounded for all $d$. In particular, we can use $M=5$ for $d\ge 41$ and $M=3$ for $d\ge 316$, and for all but very low values of $d$ the term $\log_2((4d+4)/\delta)+M+1$ is dominated by the term $\log_2(4d+4)$.

\begin{proof}
For each root $\alpha_i$, there is a point $z^{(i)}\in\Sd$ in the central subchannel of the largest channel of $\alpha_i$, and its orbit is $R$-central and satisfies $d_U(z^{(i)},N_p(z^{(i)}))\le\tau<2\log d$ (Theorem~\ref{Thm:EfficientGrid}). The claim thus follows from Proposition~\ref{Prop:IsolatedRoots} and the remark thereafter.
\end{proof}

\begin{remark}[Expected mutual distance between roots]
The results in this section were under the assumption that all (or at least some) roots were $\delta$-separated for some $\delta>0$. If $d$ roots are placed independently and randomly into $\disk$ (with respect to planar Lebesgue measure), then the mutual distance between any two roots is easily seen to be at least $O(1/d)$. Theorem~\ref{Thm:EfficientGrid} thus applies and yields, for fixed $\delta$, a number of iterations of at most $O(d^3\log^3 d+d\log|\log\eps|)$.

If not the locations of the roots are chosen randomly, but for instance the coefficients, then the roots may no longer be equidistributed with respect to area; they tend to distribute uniformly along a circle \cite{ErdoesTuran}, and the expected mutual distance is at least $O(1/d^{2})$. In any case, the relation between coefficients and roots is algebraic, so the expected mutual distance $\delta$ between roots is bounded by a power law in $d$, say $\delta\ge 1/d^\beta$ with some $\beta\ge 1$, but since our estimates only involve $\log|\delta|$, this still becomes only a constant factor $\beta$ in the number of iterations.
\end{remark}

\begin{remark}[Further improvements]
The greatest loss in our estimates is in the most basic of our estimates, in Lemma~\ref{Lem:LogDerivative}: if $|z-N_p(z)|<s$, then $|z-\alpha|<ds$ for some root $\alpha$. This bound is sharp only if all roots form a single multiple root, and then indeed the distance to the root is multiplied by $(d-1)/d$ in each Newton iteration. If the roots are randomly distributed, then the bound is much better, and this leads to significant improvements. Refining our methods in this direction, the following is shown in \cite{NewtonTodor}: \emph{If the $d$ roots are distributed independently in $\disk$ and randomly with respect to Lebesgue measure of $\disk$, then the number of iterations for the same grid $\Sd$ as before is at most $O(d^2\log^4d+d\log|\log\eps|)$, with high probability.} This improves our bound by a factor of $d/\log d$ and it is optimal except for some powers of $\log$: if we have $d$ starting points outside of $\disk$ at radius $r_0>e$, then the simple estimate $((d-1)/d)^d\approx 1/e$ implies that each of them takes approximately $d$ iterations to move from any radius $r$ to radius $r/e$, so the $d$ points together need $O(d^2)$ iterations even to get close to $\disk$ (and if $r_0=e^\beta$ with $\beta\in[0,1]$, then only a constant factor $\beta$ is gained). Our universal set of starting points requires us to place the starting points uniformly outside of $\disk$, and under this assumption the number of iterations is essentially best possible. 
This remark also applies when the roots are randomly distributed along a circle, for instance when the coefficients are independently randomly distributed.

As mentioned in the introduction, the expectations of efficient root finding are verified in practice \cite{NewtonRobin}: all $d$ roots of various polynomials of degrees up to $2^{20}>10^6$ have been found in between $3d^2$ and $6d^2$ iterations by a program based on the theory described here. 
\end{remark}

%\newpage

\section{Non-Uniformly Separated Roots}
\label{Sec:WorstCase}

If the roots are not uniformly $\delta$-separated, then they may be multiple, and the local rate of convergence may be linear rather than quadratic. For practical purposes, multiple roots are the same as simple roots at a distance smaller than the required precision $\eps$. Our previous estimate on the required number of iterations scales with $\delta$ essentially as $O(d^3|\log\delta|)$: this is of course unbounded, but diverges slowly as $\delta\to 0$; for random distributions of roots (for instance with respect to Lebesgue measure of the plane) the expected value is finite and of moderate size.

However, there are of course important polynomials with multiple or near-multiple roots. Thus we will now provide a uniform bound on the required number of iterations for all polynomials in $\Pd$. We will assume that all roots are simple, but since we do not assume a lower bound on their mutual distance, the estimates hold for multiple roots as well, by continuity.
The point grid that we will use is the same as before. 

The issue of ``clusters of roots'' is relevant from many points of view: from a distance, such clusters look like multiple roots (resulting in slowing down the Newton dynamics), and only near such a cluster does the dynamics begin to see the roots separately (in fact, sufficiently far outside any disk containing all roots of a degree $d$ polynomial, the roots look like a single root of multiplicity $d$, which explains the linear convergence with the factor $(d-1)/d$ near $\infty$). For a systematic study of clusters of roots, as well as further references on this topic, see \cite{GLSY} (however, in this reference the assumption is made that the number of roots within any cluster is known ahead of time). 

Our estimates will be based on the following stopping criterion.

\begin{lemma}[Worst-case stopping criterion] 
\label{Lem:StoppingWorstCase} \lineclear
If $z\in U_\alpha$, the immediate basin of a root $\alpha$, 
then $|z-\alpha|<f_d |z-N_p(z)|$, where 
\[
f_d=\frac{d^2(d-1)}{2(2d-1)}{{2d}\choose{d}} < d^2 4^{d-1}
\]
depends only on $d$ and satisfies 
$\log_2 f_d<2(d-1)+ 2\log_2 d$.
\end{lemma}
\begin{proof}
This result is proved in \cite[Lemma~5]{NewtonIterations}, using an iterated ``cluster of roots'' argument. 
\end{proof}

The difficulty in this result is the following: if $|z-N_p(z)|<\eps/d$, then $z$ is $\eps$-close to some root $\alpha'$ by Lemma~\ref{Lem:LogDerivative}; but even if $z\in U_\alpha$, this does not mean that $\alpha'=\alpha$. Stopping the iteration at $z$ and declaring $z$ as an approximation to a nearby root, which is necessarily $\alpha'$, runs the danger that $z$ was the only root guaranteed to find $\alpha$, and the algorithm might miss $\alpha$ altogether. The iterated cluster of roots argument in \cite{NewtonIterations} argues that either $\alpha'=\alpha$, or some other roots $\alpha''$ must be reasonably close to $z$ (with constants depending on $d$). Then either $\alpha''=\alpha$, or a further root $\alpha'''$ must be close, etc. In order to assure that $z$ is indeed $\eps$-close to the ``correct'' root $\alpha$, in the (unlikely) iterated worst case we are led to the factor $f_d$ in Lemma~\ref{Lem:StoppingWorstCase}.

The following result is the ``worst-case'' version of Proposition~\ref{Prop:IsolatedRoots}. 

%\newpage

\begin{proposition}[Number of iterations on orbit, worst case]
\label{Prop:NumberIterationsWorstCase} \lineclear
If $(z_n)$ is an $R$-central orbit in $U$ with $d_U(z_0,z_{1})\le \tau$ and $|z_0|\le R\le 2^M-1$, then $|z_N-\alpha|<\eps$ provided $N$ is at least
\begin{align*}
&\pi (4d+1)^2(2\tau d+1/4)\lceil 2\tau+4\log 3\rceil \left(2(d-1)+\log_2d-\log_2(\eps)+M+1\rule{0pt}{11pt}\right) \\
& \quad +\log_2|\log_2\eps-5|+1
\\
& = \left(128 \pi d^4\tau^2 +64\pi d^3
\tau^2 |\log_2\eps|  \right) (1+o(1))+\log_2|\log_2\eps| 
\\
& \in O\left(d^4\tau^2+d^3\,\tau^2(|\log\eps|)+\log|\log\eps|\right)
\;.
\end{align*}
\end{proposition}
\begin{proof}
We iterate the orbit $(z_n)$ while $|z_n-z_{n+1}|\ge \eps/f_d$ and  $|\phi(z_n)|\ge e^{1/2}-1$. If at some time $|z_n-z_{n+1}|< \eps/f_d$, then we can stop by Lemma~\ref{Lem:StoppingWorstCase}, and if $|\phi(z_N)|< e^{1/2}-1$, then Lemma~\ref{Lem:HypDistanceOrbit} (2) applies and only $\log_2|\log_2\eps-5|$ further iterates are required until $\eps$-precision is reached.

We now estimate how many iterates are necessary until $|z_N-z_{N+1}|< \eps/f_d$; we may suppose that along the way, we always have $|\phi(z_n)|>e^{1/2}-1$. 

\looseness -1
We will use Corollary~\ref{Cor:NumberIteratesS_k} for $k=-M,0,1,2,\dots,K$, where $K$ is the least integer such that $1/d2^{K+1}\le\eps/f_d$, i.e., $K= \lceil  
\log_2(f_d/2d\eps)\rceil <2(d-1)+\log_2d-\log_2(2\eps)+1$.
The number of orbit points in $S_k$ is at most 
$\pi (4d+1)^2(2\tau d+1/4)\lceil 2\tau+4\log 3\rceil $.
For the final value $k=K$, the total number of points in $S^*_K$ with $|z_n-z_{n+1}|\ge\eps/f_d\ge 1/d2^{K+1}$ satisfies the same bound (Proposition~\ref{Prop:NumberOrbitPoints}). 
Moreover, we have $D_{2^M-1}(0)\subset S_{-M}\cup S_{-M+1}\cup \dots\cup S_{K-2}\cup \dots\cup S_{K-1}\cup S^*_K$, so the total number of iterations with $|z_n-z_{n+1}|>\eps/f_d$ is at most  $K+M+1$ times the number for each $S_k$, and by hypothesis the orbit never leaves the disk $D_R(0)\subset D_{2^M-1}(0)$.
\end{proof}

\goodbreak

\begin{theorem}[Worst case number of iterations]
\lineclear 
For each degree $d$, the set $\Sd$ constructed in Definition~\ref{Def:EfficientGrid} has the following universal property. If $p$ is any complex polynomial, normalized so that all its roots are in $\disk$, and with simple or multiple roots at arbitrary mutual distances, then there are $d$ points in $\Sd$ that converge to the $d$ roots of $p$, and so that the combined number of iterations required to reach $\eps$-precision is at most 
\begin{align*}
&\pi (4d+1)^2(2\tau d+1/4)\lceil 2\tau+4\log 3\rceil \left(2(d-1)+\log_2d-\log_2(\eps)+M+1\rule{0pt}{11pt}\right) \\
& \quad +d\log_2|\log\eps-5|+1
\\
&
= \left(512\pi d^4\log^2d + 256\pi d^3\log^2d|\log\eps| \right) (1+o(1))
\\
& \in O\left(d^4\log^2 d+d^3\,\log^2d|\log\eps|\right)
\;.
\end{align*}
\end{theorem}
\begin{proof}
It suffices to prove this for the case that all roots are simple; the case of multiple roots follows by continuity. 
The set $\Sd$ intersects the central subannulus of the channel with largest modulus of each root in at least one point, so the corresponding orbits are $R$-central and all their orbit points $z$ satisfy $d(z,N_p(z))<\tau:=2\log d$. 

By Proposition~\ref{Prop:NumberIterationsWorstCase}, each of these points needs at most
\begin{eqnarray*}
&&\pi (4d+1)^2(2\tau d+1/4)\lceil 2\tau+4\log 3\rceil \left(2(d-1)+\log_2d-\log_2(2\eps)+M+2\rule{0pt}{11pt}\right) \\
&& \quad +\log_2|\log_2\eps-5|+1
\end{eqnarray*}
iterations to be $\eps$-close to the corresponding root. Since all roots again have to compete for the area within $\disk$, the total number of iterations combined to get $\eps$-close to all $d$ roots, for one starting point per root, satisfies the same bound, except that a factor $d$ comes in in the part of the estimate where the roots do not compete for area, and this is the term with $\log_2|\log_2\eps|$ (which is subordinate to the $d\,|\log\eps|$ term). This proves the claim.
\end{proof}

%\newpage

\begin{remark}
We believe that this result in the worst case can be improved at least by a factor of $d$: the factor $f_d$ in Lemma~\ref{Lem:StoppingWorstCase} is exponential in $d$, and the worst case leading to this estimate seems very unrealistic. Even though $f_d$ enters only logarithmically, $\log f_d$ still contributes a factor of $d$. If $f_d$ could be replaced by a polynomial in $d$, this would gain a factor of $d/\log d$. (However, the complexity in $\eps$ really is $|\log\eps|$, rather than $\log|\log\eps|$, in the presence of multiple roots because Newton's method at multiple roots converges only linearly, not quadratically). 
\end{remark}

%\newpage

\section{Some Numerical Experiments}
\label{Sec:Experiments}

In this section we briefly report on some numerical experiments performed jointly with Robin Stoll. We mentioned earlier one set of experiments, strictly based on the theory described here, that manages to find all roots of various polynomials of degree up to a million with complexity between $3d^2$ and $6d^2$, confirming our estimates in practice. 

All these are worst case bounds,  and the optimality of the results in this manuscript is based on these worst case assumptions. Here we present the results of a ``more optimistic'' implementation of Newton's method. 
We investigate two families of polynomials of degrees $2^n$ with $n\le 27$, i.e.\ of degrees up to $134$ million. These polynomials were chosen so that they and their derivatives can be evaluated efficiently by recursion (our focus is on root finding, not on polynomial evaluations): for a given quadratic polynomial $q$, we find the periodic points of period $n$; these are roots of the polynomial $q^{\circ n}(z)-z$ (where $q^{\circ n}$ denotes again the $n$-th iterate). Our two families of polynomials correspond to periodic points of $q_2(z)=z^2+2$ and $q_i(z)=z^2+i$.

\begin{table}[htbp]
\textbf{Periodic points of $z^2+2$.} 
{\small \framebox{
\begin{tabular}{rrrr | rl}
period & degree & iterations &  iterations & computing  \\
& & & $d\log^2d$ & time (secs) \\
\hline 
12 &	4\,096		&$756\ d$	 	&2.52 & 1	\\
13 &	8\,192		&$798\ d$		&2.27 & 3 	\\
14 &	16\,384		&$1053\ d$	&2.58 & 7	\\
15 &	32\,768		&$1220\ d$	&2.61 & 16	\\
16 &	65\,536		&$1399\ d$	&2.63 & 37	\\
17 &	131\,072		&$1585\ d$	&2.64 & 87	\\
18 &	262\,144		&$1786\ d$	&2.65 & 201	\\
19 &	524\,288		&$1988\ d$	&2.65 & 462	\\
20 &	1\,048\,576	&$2210\ d$	&2.65 & 1\,058	\\
21 &	2\,097\,152	&$2437\ d$	&2.66 & 2\,407	\\
22 &	4\,194\,304	&$2678\ d$	&2.66 & 5\,453  \\
23 &	8\,388\,608	&$2945\ d$	&2.67 & 18\,520 &  \\
24 & 16\,777\,216	& $3204\ d$	&2.67 & 32\,401 &	\\
25 & 33\,554\,432 	& $3457\ d$	&2.66 & 34\,500 &	\\
26 & 67\,108\,864	& $3738\ d$	&2.66 & 76\,698 & 	\\
27 & 134\,217\,728	& $4044\ d$	&2.67 & 320\,567 & (89 hours 3 min) 	 \\
\end{tabular}}
\medskip
\caption{Finding periodic points of $q_2(z)=z^2+2$. The first  columns show period $n$ and degree $d=2^n$ of $q_2^{\circ n}(z)-z$. The third column shows the number of Newton iterations required to find all $d$ roots.  The next column shows that the number of iterations seems to converge to $2.67\ d\log^2d$. The final column shows the computing time on a standard PC computer from about 2012 (single core): this time seems to scale with $d\log^2d$ as well (with some variations). }
\label{Tab:Experiment_2}
}
%\end{table}

%\begin{table}[htbp]
\textbf{Periodic points of $q_i(z)=z^2+i$.} 
{\small
\framebox{
\begin{tabular}{rrrr  r | rl l}
 period &  degree &  iterations &  iterations &  iterations & computing & \!time \\
& & & $d\log^2d$  & $d\log d$ & (secs) \ \ \ \  \\
\hline 
12 &	4\,096		& $318\ d$ & 1.06 & 26.50 &  1	\\
13 &	8\,192		& $351\ d$ & 1.00 & 27.00  & 1	\\
14 &	16\,384		& $385\ d$ & 0.94 & 27.50  &  3	\\
15 &	32\,768		& $418\ d$ & 0.89 & 27.87  & 7	\\
16 &	65\,536		& $451\ d$ & 0.85 & 28.19  & 15	\\
17 &	131\,072		& $485\ d$ & 0.81& 28.53  &  33	\\
18 &	262\,144		& $518\ d$ & 0.77 & 28.78  &  71	\\
19 &	524\,288		& $551\ d$ & 0.73 & 29.00  & 153	\\
20 &	1\,048\,576	& $585\ d$ & 0.70 & 29.25 & 332 	\\
21 &	2\,097\,152	& $618\ d$ & 0.67& 29.43 & 716    	\\
22 &	4\,194\,304	& $652\ d$ & 0.65& 29.64 & 1\,541 \\
23 &	8\,388\,608	& $685\ d$	&0.62 & 29.78  & 3\,309 \\
24 & 16\,777\,216	& $718\ d$	&0.60 & 29.92  & 7\,091	\\
25 & 33\,554\,432 	& $752\ d$	&0.58 & 30.08  & 15\,139	\\
26 & 67\,108\,864	& $785\ d$	&0.56 & 30.19  & 32\,325	\\
27 & 134\,217\,728	& $818\ d$	&0.54 & 30.30   & 69\,302 & (19 h 15 min)\\
\end{tabular}}
\medskip
\caption{The same experiment for periodic points of $q_i(z)=z^2+i$. This time, the number of Newton iterations seems to scale with $d\log d$ for degrees up to and exceeding 134 million. } 
\label{Tab:Experiment_i}
}
\end{table}

It turns out that for both families (and numerous others), it is quite possible to find all roots for degrees many millions, and with a guarantee that all roots were indeed found. No issues about large computing precision were encountered (all our polynomials have simple roots). Moreover, the required number of Newton iterations scales with $d\log^2 d$ or even $d\log d$, and roughly the same applies to the computing time (on a PC computer of about 2012). The overall outcome is listed in Tables~\ref{Tab:Experiment_2} and \ref{Tab:Experiment_i}.
The details of these and further experiments are reported in \cite{NewtonRobin,NewtonRobin2}. We show them here in order to support our claim that Newton's method is a root finding method that has both good theory (as developed in this text and its subsequent improvements in \cite{NewtonTodor}) and that works remarkably well in practice. 

In particular, for both families of polynomials all roots could verifiably be found on a standard PC for degree $134$ million in hours, or a few days,  of computing time (on a single core computer; multiple cores would reduce the time almost linearly). Even for these extremely large degrees, no special arithmetic was required, nor did the software have to be adapted. It seems that the limiting factor for the computations are still not numerical issues but RAM memory constraints (for degree $2^{27}$ the output file alone required about 9 Gigabytes of memory); quite possibly the limit can be pushed substantially further by optimizing the software. 

It is a noteworthy and perhaps amusing fact that for $q_2^{\circ n}(z)-z$ the constant coefficient has magnitude greater than $2^{2^n}$, far greater than can be stored using any standard arithmetic (for $n=27$ this quantity has more than $40$ million decimal digits!, while all roots are clustered within a disk of radius $2$ and are thus very close to each other) --- but this caused no problem for the computations at all (our iterative evaluation of the polynomials does not require the coefficients). 

%\newpage

\appendix
\section{Geometry of Hyperbolic Geodesics}
\label{Sec:GeometryHypGeodesics}

In this appendix, we will prove the claim that the line segments $X(t)$ as introduced in Section~\ref{Sec:AreaPerIteration} are disjoint.
We will repeatedly use Ahlfors' theorem that every Euclidean disk $D\subset U$ is convex with respect to hyperbolic geodesics in $U$  \cite{Jorgensen}. Hence for every geodesic $\gamma$ the set $\gamma\cap \ovl D$ is connected.

\begin{lemma}[Euclidean curvature radius]
\label{Lem:HypCurvatureRadius}\lineclear
For every $t_0\in\R$, the Euclidean curvature radius of $\gamma$ at $\gamma(t_0)$ is at least $\eta(t_0)/2$.
\end{lemma}
\begin{proof}
Let $R\in(0,\infty)$ be the Euclidean curvature radius of $\gamma$ at $\gamma(t_0)$ (if $R=\infty$, then we have nothing to show). If $R<\eta(t_0)/2$, then let $C$ be a circle of radius between $R$ and $\eta(t_0)/2$ and tangent to $\gamma$ at $\gamma(t_0)$, and so that for $t\neq t_0$ sufficiently close to $t$, the point $\gamma(t)$ is in the disk bounded by $C$; call this disk $D$. But then $\gamma(t_0)$ disconnects $\gamma\cap U$ in contradiction to Ahlfors' theorem that disks are hyperbolically convex.
\end{proof}

\begin{remark}
This bound might well be sharp. (Bj\"orn Gustafsson \cite[Corollary~8.6]{CurvatureBound}  observed that it is sharp for domains $U\subset\Cbar$ that may contain the point at $\infty$, and Edward Crane observed that it is not far from being sharp for domains $U=\C\sm\R^-_0)$.
\end{remark}

Recall that for a Riemann domain $U\subset\C$ and a hyperbolic geodesic $\gamma\colon\R\to U$ parametrized by Euclidean arc length, we defined $\eta(t)$ as the Euclidean distance of $\gamma(t)$ to $\partial U$, and $X(t)$ as the straight line segment (without endpoints) of length $\eta(t)$ with center at $\gamma(t)$ that intersects $\gamma$ at $\gamma(t)$ in a right angle.

\begin{proposition} %[Disjoint $X(t)$]
\label{Prop:All_Xt_Disjoint} %\lineclear
All $X(t)$ are disjoint.
\end{proposition}
\begin{proof}
(0) Suppose $X(t_0)$ and $X(t_1)$ intersect; without loss of generality, suppose that $\eta(t_0)\ge\eta(t_1)$ and $t_1>t_0$. Let $D_0$ be the open disk centered at $\gamma(t_0)$ and with radius $\eta(t_0)$, and let $C_1$ and $C_2$ be the two circles of radius $\eta(t_0)/2$ tangent to $\gamma'(t_0)$; then both circles are tangent (from the inside) to $\partial D_0$, and $X(t_0)$ is exactly the open straight line segment connecting their centers. Let $D_1$ and $D_2$ be the two open disks bounded by $C_1$ and $C_2$, and let $\ell$ be the straight line through their centers. Without loss of generality, we may assume that $X(t_0)$ and $X(t_1)$ intersect within $D_1$.

\begin{figure}[h]
\includegraphics[scale=0.7]{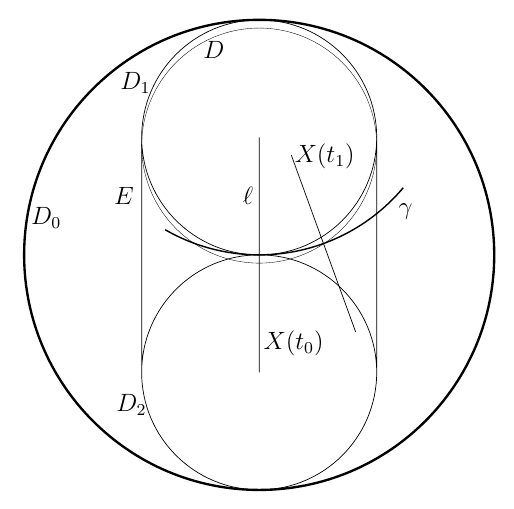}
\caption{Illustration of the proof of Proposition~\ref{Prop:All_Xt_Disjoint}. The two disks $D_1$ and $D_2$ are tangent to $\gamma(t)$ at $t=t_0$, the set $E$ is the convex hull of $D_1$ and $D_2$. Within $E$, the curve $\gamma$ is a single graph over the horizontal direction. The disk $D$ (marked by a thin boundary circle) is a vertical translate of $D_1$ so that $\gamma(t_1)\in\partial D$. The fact that $\gamma$ must leave $D$ at $\gamma(t_1)$ forces the slope of $\gamma$ at $t_1$ to be no larger than the slope of $\partial D$ at that point, and this tilts $X(t_1)$ so that its distance to $\ell$ is greater than the radius of $D_1$. }
\end{figure}

(1) First observe that the geodesic $\gamma$ is disjoint from $D_1$ and $D_2$: if there is some point $\gamma(t_2)\in D_1$, say, then consider the disk $D'_1\subset D_1$ so that $\partial D'_1\cap\partial D_1=\{\gamma(t_0)\}$ and $\gamma(t_2)\in\partial D'_1$. Then $\ovl D'_1\subset U$ and the radius of $D'_1$ is smaller than the radius of $D_1$, and hence smaller than the radius of curvature of $\gamma$ at $\gamma(t_0)$ by Lemma~\ref{Lem:HypCurvatureRadius}, so $\gamma(t)$ is not in $\ovl D'_1$ for $t$ near $t_0$, except that $\gamma(t_0)\in\partial D'_1$. But $\gamma(t_2)\in\partial D'_1$, and this contradicts Ahlfors' hyperbolic disks theorem.

(2) Let $E$ be the convex hull of $D_1\cup D_2$. We claim that $\gamma\cap \ovl E$ is connected. Indeed, suppose $\gamma(t_2)\not\in \ovl E$, but $\gamma(t_3)\in \ovl E$, with $t_0<t_2<t_3$, say. Then there is a disk $D'$ with radius $\eta(t_0)/2$ and with center on $X(t_0)$, and so that $\gamma(t_3)\in \ovl D'$; it satisfies $\ovl D'\subset U$. But since $\gamma(t_0)\in\ovl D'$, it follows that $\gamma\cap\ovl D'$ is not connected, a contradiction. 

(3) Our next claim is that along $\gamma\cap E$, no tangent vector is  perpendicular to $\gamma'(t_0)$; this means that, if we call the direction of the tangent vector $\gamma'(t_0)$ ``horizontal'', then $\gamma\cap E$ is a graph over the horizontal direction. For a proof by contradiction, suppose there is a $t_2>t_0$, say, with $\gamma(t_2)\in E$ and $\gamma'(t_2)$ perpendicular to $\gamma'(t_0)$. Let $D'_1$ and $D'_2$ be the two disks with centers on $\ell$ and with radii $\eta(t_0)/2$ and so that $\gamma(t_2)\in\partial D'_1\cap\partial D'_2$. Since $\gamma(t_2)\not\in D_1\cup D_2$, it follows that $\gamma(t_0)\in D'_1\cap D'_2$, and hence that $D'_i\subset E$ for $i=1,2$. One of the two disks $D'_1$ and $D'_2$ thus has the property that $\gamma(t_2+\eps)\in D'_i$ for small $\eps>0$, but not for small $\eps<0$. Since $\gamma(t_0)\in D'_i$ and $t_0<t_2-\eps<t_2$, this contradicts Ahlfors' theorem once again.

(4) Now consider the point $\gamma(t_1)$ and let $D$ be the unique disk of radius $\eta(t_0)/2$ with center on $\ell$ and so that $\gamma(t_1)\in\partial D$ (this leaves two choices, and we take the disk with center closest to the center of $D_1$). Since $\gamma(t_1)\not\in \ovl D_1$, the center of $D$ is strictly between the centers of $D_1$ and $D_2$, and hence $\gamma(t_0)\in D$. By Ahlfors' theorem again, we have $\gamma([t_0,t_1])\subset \ovl D$, and $\gamma(t)\not\in \ovl D$ for $t>t_1$. 

Since $\gamma$ must leave $D$ in the direction of increasing $t$, the tangent vector $\gamma'(t_1)$ must either be parallel to the tangent vector of $\partial D$ at $\gamma(t_1)$, or its slope must be smaller. If they are parallel, then the distance from $\gamma(t_1)$ to $\ell$ along $X(t_1)$ is exactly $\eta(t_0)/2$ (the radius of $D$), while the length of $X(t_1)$ is $\eta(t_1)/2\le\eta(t_0)/2$ (from the center point $\gamma(t_1)$ in both directions), so $X(t_1)$ cannot intersect $X(t_0)$. If the slope of $\gamma'(t_1)$ is smaller than that of the circle, then the distance from $\gamma(t_1)$ to $\ell$ in the direction of $X(t_1)$ is even greater. Therefore, $X(t_0)$ and $X(t_1)$ cannot intersect.
\end{proof}


\begin{thebibliography}{MM}
\bibitem[AHU]{HopcroftUllman} Alfred V. Aho, John E. Hopcroft, and Jeffrey D. Ullman, \emph{The design and analysis of computer algorithms}, Addison Wesley, 1974.

\bibitem[BAS]{NewtonTodor} Todor Bilarev, Magnus Aspenberg, and Dierk Schleicher, \emph{On the speed of convergence of Newton's method for complex polynomials}. Mathematics of Computation \textbf{85} 298 (2016), 693--705. 

\bibitem[BR]{BiniRobol} Dario A. Bini and Leonardo Robol, 
\emph{Solving secular and polynomial equations: a multiprecision algorithm}. Journal of Computational and Applied Mathematics \textbf{272} 276--292 (2014).

\bibitem[BBEGG]{BiniEigenvalue}
Dario A. Bini, Paula Boito, Yuli Eidelman, Luca Gemignani, and Israel Gohberg, 
\emph{A fast implicit $QR$ eigenvalue algorithm for companion matrices}.
Linear Algebra Appl. \textbf{432} 8 (2010), 2006--2031. 

\bibitem[BLS]{BLS} B\'ela Bollob\'as, Malte Lackmann, and Dierk Schleicher, \emph{A  small probabilistic universal set of starting points for finding roots of complex polynomials by Newton's method.} Mathematics of Computation, \textbf{82} 281 (2013), 443--457.


\bibitem[ET]{ErdoesTuran}
Paul Erd\H{o}s and P\'al Tur\'an, 
\emph{On the distribution of roots of polynomials}.
Ann. of Math. (2) \textbf{51} (1950), 105--119. 

\bibitem[GLSY]{GLSY}
Marc Giusti, Gr\'egoire Lecerf, Bruno Salvy, and Jean-Claude Yakoubsohn, \emph{Location and approximation of clusters of zeros of analytic functions}. Foundations of Computational Mathematics \textbf{5} (2005), 257--311.


\bibitem[GS]{CurvatureBound}%\looseness-1
Bj\"orn Gustafsson and Ahmet Sebbar,
\emph{Critical points of Green's function and geome\-tric function theory}. 
Indiana Univ. Math.~J.~\textbf{61} 3, 939--1017 (2012). 	arXiv: 0912.1223.



\bibitem[HSS]{NewtonInventiones} \looseness-1
John Hubbard, Dierk Schleicher, and Scott Sutherland, {\em How to find all roots of complex polynomials by Newton's method}. Inventiones Mathematicae {\bf 146} (2001), 1--33.

\bibitem[J]{Jorgensen}
Vilhelm J{\o}rgensen, \emph{On an inequality for the hyperbolic measure and its applications in the theory of functions}. Math. 
Scand. \textbf 4 (1956) 113--124. 

\bibitem[MS]{NewtonSebastian}
Sebastian Mayer and Dierk Schleicher, {\em Immediate and virtual basins of Newton's method for entire functions}. Annales Institut Fourier, Grenoble {\bf 56} 2 (2006), 325--336.

\bibitem[MMS]{NewtonAlgorithm}
Khudoyor Mamayusupov, Sabyasachi Mukherjee, and Dierk Schleicher, 
\emph{Turning Newton's Method into an Algorithm with Predictable Complexity}.
Manuscript, in preparation.


\bibitem[McN1]{McNamee} John M. McNamee, \emph{A 2002 update of the supplementary bibliography on roots of polynomials}. J. Comput. Appl. Math. \textbf{142} 2 (2002), 433--434.

\bibitem[McN2]{McNameeBook}
John M. McNamee, \emph{Numerical methods for roots of polynomials. Part I. Studies in Computational Mathematics} \textbf{14}. Elsevier B. V., Amsterdam, 2007.

\bibitem[MP]{McNameePanBook}\looseness-1
John M. McNamee and Victor Pan, \emph{Numerical methods for roots of polynomials. Part II. Studies in Computational Mathematics} \textbf{16}. Elsevier B. V., Amsterdam, 2013.

\bibitem[MR]{HelenaLasse} Helena Mihaljevi\'c-Brandt and Lasse Rempe-Gillen,
\emph{Absence of wandering domains for some real entire functions with bounded singular sets}. Mathematische Annalen \textbf{357} 4 (2013), 1577--1604.

\bibitem[MB]{MoenckBorodin}
Robert T.\ Moenck and Allan B.\ Borodin,
\emph{Fast Modular Transform via Division},
{Proc. 13th annual symposium on switching and automata theory} 
90--96, IEEE Comp. Society Press, Washington, DC, 1972.


\bibitem[P1]{Pan}
Victor Pan, 
\emph{Approximating Complex Polynomial Zeros: Modified Weyl's
Quadtree Construction and Improved Newton's Iteration}. 
Journal of Complexity \textbf{16} (2000), 213--264.

\bibitem[P2]{Pan_NearOptimal}
Victor Pan, 
\emph{Univariate Polynomials: Nearly Optimal Algorithms for Numerical Factorization and Root-finding}. 
J. Symbolic Computation \textbf{33} 5 (2002), 701--733.


\bibitem[Pr]{Prz}
Feliks Przytycki, {\em Remarks on the simple connectedness of basins of sinks for
iterations of rational maps}. In: Dynamical Systems and Ergodic Theory, ed. by
K. Krzyzewski, Polish Scientific Publishers, Warszawa (1989), 229--235.

\bibitem[R]{Renegar}
James Renegar,
\emph{On the worst-case arithmetic complexity of approximating zeros of polynomials}.
J. Complexity \textbf{3} 2 (1987), 90--113. 


\bibitem[R\"u]{JohannesSurvey}
Johannes R\"uckert, {\em Rational and transcendental Newton maps}. 
In: \emph{Holomorphic dynamics and renormalization. A volume in honour of John Milnor's 75th birthday}, ed. Mikhail Lyubich and Michael Yampolsky, Fields Institute Communications \textbf{53} (2008), 197--212. 

\bibitem[RS]{NewtonJohannes} 
Johannes R\"uckert and Dierk Schleicher, {\em On Newton's method for entire functions}. Journal of the London Mathematical Society {\bf 75} 3 (2007), 659--676. 

\bibitem[Sch1]{NewtonIterations} 
Dierk Schleicher, {\em On the number of iterations of Newton's method for complex polynomials}. Ergodic Theory and Dynamical Systems {\bf 22} (2002), 935--945.

\bibitem[Sch2]{NewtonFields}
Dierk Schleicher, \emph{Newton's
method as a dynamical system: efficient root finding of polynomials
and the Riemann $\zeta$ function}. In: \emph{Holomorphic dynamics and renormalization. A volume in honour of John Milnor's 75th birthday}, ed. Mikhail Lyubich and Michael Yampolsky, Fields Institute Communications \textbf{53} (2008), 213--224.

\bibitem[SSt]{NewtonRobin}
Dierk Schleicher and Robin Stoll,
\emph{Newton's method in practice: finding all roots of polynomials of degree one million efficiently}. Journal of Theoretical Computer Science, to appear (2015).

\bibitem[SSt2]{NewtonRobin2}
Dierk Schleicher and Robin Stoll,
\emph{Newton's method in practice II: near-optimal complexity for 
finding all roots of some polynomials of large degrees}.
Manuscript, in preparation.

\bibitem[Sh]{Mitsu}
Mitsuhiro Shishikura, 
{\em The connectivity of the Julia Set and fixed points}. In: \emph{Complex dynamics: families and friends}, ed. Dierk Schleicher, AK Peters, Wellesley/MA, 2009, 257--276.

\bibitem[Sm]{Smale}
Steve Smale, {\em Newton's method estimates from data at one point}. In: The merging of disciplines:
new directions in pure, applied, and computational mathematics (Laramie, Wyo.,
1985), 185--196, Springer, New York, 1986.

\end{thebibliography}
\end{document}